\documentclass[sn-basic,iicol]{sn-jnl}% Basic Springer Nature Reference Style/Chemistry Reference Style
\usepackage{mathtools}
\usepackage{graphicx}%
\usepackage{multirow}%
\usepackage{amsmath,amssymb,amsfonts}%
\usepackage{amsthm}%
\usepackage{mathrsfs}%
\usepackage{rotating}
\usepackage{etoolbox}
\usepackage{url}
\usepackage[title]{appendix}%
\usepackage{xcolor}%
\usepackage{textcomp}%
\usepackage{manyfoot}%
\usepackage{booktabs}%
\usepackage{algorithm}%
\usepackage{algorithmic}%
\usepackage{listings}%
\usepackage{geometry}
\usepackage{hyperref}
\usepackage{breakurl}
\usepackage{caption}
\usepackage{subcaption}
\usepackage{wrapfig}
\makeatletter
\patchcmd{\NAT@citex}
  {\@citea\NAT@hyper@{%
     \NAT@nmfmt{\NAT@nm}%
     \hyper@natlinkbreak{\NAT@spacechar\NAT@@open\if*#1*\else#1\NAT@spacechar\fi}%
       {\@citeb\@extra@b@citeb}%
     \NAT@date}}
  {\@citea\NAT@nmfmt{\NAT@nm}%
   \NAT@spacechar\NAT@@open\if*#1*\else#1\NAT@spacechar\fi\NAT@hyper@{\NAT@date}}
  {}{}
  \patchcmd{\NAT@citex}
  {\@citea\NAT@hyper@{%
     \NAT@nmfmt{\NAT@nm}%
     \hyper@natlinkbreak{\NAT@aysep\NAT@spacechar}{\@citeb\@extra@b@citeb}%
     \NAT@date}}
  {\@citea\NAT@nmfmt{\NAT@nm}%
   \NAT@aysep\NAT@spacechar\NAT@hyper@{\NAT@date}}{}{}
\makeatother  
\newtheoremstyle{plain}
  {8pt}   % ABOVESPACE
  {8pt}   % BELOWSPACE
  {\itshape}  % BODYFONT
  {0pt}       % INDENT (empty value is the same as 0pt)
  {\bfseries} % HEADFONT
  {.}         % HEADPUNCT
  {5pt plus 1pt minus 1pt} % HEADSPACE
  {}          % CUSTOM-HEAD-SPEC
\theoremstyle{plain}
\newtheorem{theorem}{Theorem}
\newtheorem{proposition}[theorem]{Proposition}
\theoremstyle{definition}
\newtheorem{remark}{Remark}%
\theoremstyle{plain}
\newtheorem{definition}{Definition}%

\DeclareMathOperator*{\argmin}{arg\,min}
\begin{document}

\title[Article Title]{Uncertainty Quantification for Linear Inverse Problems with Besov Prior: A Randomize-Then-Optimize Method}


\author*[1]{\fnm{Andreas} \sur{Horst}}\email{ahor@dtu.dk}

\author[2]{\fnm{Babak} \sur{Maboudi Afkham}}\email{babak.maboudi@oulu.fi}

\author[1]{\fnm{Yiqiu} \sur{Dong}}\email{yido@dtu.dk}
\author[1]{\fnm{Jakob} \sur{Lemvig}}\email{jakle@dtu.dk}

\affil[1]{\orgdiv{Department of Applied Mathematics and Computer Science}, \orgname{Technical University of Denmark}, \orgaddress{\street{Matematiktorvet}, \city{Kgs. Lyngby}, \postcode{2800},\state{} \country{Denmark}}}
\affil[2]{\orgdiv{Research Unit of Mathematical Sciences}, \orgname{University of Oulu}, \orgaddress{\street{Pentti Kaiteran katu 1}, \city{Linnanmaa}, \postcode{}\state{} \country{Finland}}}

\subtitle{Published in Statistics and Computing. Final version available at: \url{https://doi.org/10.1007/s11222-025-10638-2}}
\abstract{In this work, we investigate the use of Besov priors in the context of Bayesian inverse problems. The solution to Bayesian inverse problems is the posterior distribution which naturally enables us to interpret the uncertainties. Besov priors are discretization invariant and can promote sparsity in terms of wavelet coefficients. We propose the randomize-then-optimize method to draw samples from the posterior distribution with Besov priors under a general parameter setting and estimate the modes of the posterior distribution. The performance of the proposed method is studied through numerical experiments of a 1D inpainting problem, a 1D deconvolution problem, and a 2D computed tomography problem. Further, we discuss the influence of the choice of the Besov parameters and the wavelet basis in detail, and we compare the proposed method with the state-of-the-art methods. The numerical results suggest that the proposed method is an effective tool for sampling the posterior distribution equipped with general Besov priors.
}




\maketitle

\section{Introduction} \label{sec1}
Inverse problems are the processes of computing unknown parameters from indirect measurements of a system. They are typically ill-posed, which means that the stability of the solution is sensitive to model inaccuracy, limited measurements, and noise in the measurements. To gain stable solutions we use regularization methods. Typically, regularization methods penalize certain function properties, such as smoothness, piecewise structure, or sparsity with respect to a specific basis.

In this paper, we take a Bayesian statistical approach to regularization. The Bayesian approach describes the inverse problem in a probabilistic manner. This enables regularization using a prior distribution on the unknown parameter and describes the randomness in both the model and measurements through the likelihood. In the Bayesian approach, the goal is to estimate a posterior distribution, i.e., the conditional distribution of the unknown parameters given some measurements. The statistical properties of the posterior can be interpreted as the level of uncertainty in the parameter, which is a key advantage of Bayesian methods compared to classical regularization methods.

Posteriors arising from Bayesian inverse problems are usually not known in closed form and cannot be sampled directly. Consequently, indirect sampling methods have to be applied to estimate moments of the posterior. In Bayesian inverse problems, Markov chain Monte Carlo (MCMC) methods are commonly used to explore the posterior distribution. MCMC methods generate samples from an ergodic Markov chain whose stationary distribution is the target posterior. Many MCMC methods, e.g., Metropolis-Hastings \citep{hastings1970}, preconditioned-Crank-Nicolson \citep{cotter2013}, and Metropolis-adjusted Langevin algorithm \citep{girolami2011} are proposal-based methods; that is, the samples generated from the Markov chain are used as proposals and accepted with some probability depending on the previous sample in the chain. This dependence makes the Markov chain samples autocorrelated which may lead to poor estimates of the posterior. Hence, constructing effective proposal-based MCMC methods that generates proposals with high posterior probability and low autocorrelation among accepted samples is an important endeavour.

In many inverse problem applications, the unknown parameters are defined on a continuum. Naive discretization of such parameters can result in discretization-dependent \citep{Matti_Lassas_2004} numerical solutions of Bayesian inverse problems, which are not guaranteed to converge to the underlying continuous unknown parameters. Infinite-dimensional Bayesian approaches guarantee the convergence to the continuous posterior and are independent of the discretization level. Despite the great success of such methods in modeling smooth parameters, addressing non-smooth or discontinuous parameters in the context of Bayesian inverse problems is still an ongoing direction of research \citep{suuronen2023bayesian,suuronen2022cauchy,uribe2022hybrid,li2022bayesian}.

%Besov priors have recently been proposed in \citep{Siltanen2013,Dashti2017} as a promising tool to incorporate sparsity in Bayesian inverse problems. 
Besov priors were introduced in \citep{Siltanen2013,Dashti2017} as a sparsity-promoting prior well-suited for a wide range of inverse problems. The properties of Besov priors under general parameter settings were initially studied in \citep{Lassas2009} and further studied in \citep{Dashti_Stuart2012}. This analysis led to important results of Besov priors such as discretization invariance and well-posedness of the posterior. These properties guarantee that the computation of the posterior is stable and that the prior properties are preserved at any discretization level. 

Besov priors with specific parameter choices have been applied mainly in imaging application such as X-ray computed tomography \citep{Rantala2006, Vanska2009, Siltanen2013, suuronen_emzir_lasanen_sarkka_roininen_2020, Cui_2021, Sakhaee_Entezari_2015, Niinimaki2007}, image deblurring \citep{Wang2017,Bui-Thanh_Ghattas_2015,Kolehmainen_2012}, and image denoising \citep{ Abramovich1998,Leporini2001}. Other recent sparsity-promoting Bayesian approaches include, but are not limited to, those in \citep{uribe2023horseshoe,kekkonen2023random,calvetti2019hierachical}.

We propose a novel infinite-dimensional Bayesian framework for inverse problems that promotes sparsity and allows for discontinuities in the solution. We use Besov priors in the general definition, which define a generalized Gaussian distribution on Besov spaces, providing control over regularity. We demonstrate how the parameters of such priors can be tuned to promote sparsity. We also show that posterior estimates are dimension independent.

We introduce a numerical method for exploring a posterior with Besov priors, based on the randomize-then-optimize (RTO) method proposed by \citep{Bardsley2014} for efficiently drawing samples from the posterior distribution. Compared to alternative MCMC approaches, the proposed samples are independent and projected towards higher probable regions of the posterior, yielding a higher acceptance rate. In addition, the RTO method does not require extensive parameter tuning, resulting in an effective method. 

The numerical experiments in this paper suggest that the proposed method provides an effective sampling method for exploring the posterior for various inverse problems, e.g., deconvolution, inpainting, and X-ray computed tomography (CT). We demonstrate that the parameters of the Besov prior can be tuned to promote the desired sparsity. In addition, we show that such priors can also exploit smoothness by choosing the right basis for the Besov space.

While numerous studies \citep{Bui-Thanh_Ghattas_2015,Wang2017,suuronen_emzir_lasanen_sarkka_roininen_2020,Sakhaee_Entezari_2015,Siltanen2013,Kolehmainen_2012} have addressed the edge-preserving properties of Besov priors under a fixed parameter setting, our work differs by considering a general parameter setup that enables a systematic study of their versatility. Our RTO sampling methodology builds upon earlier work in \citep{Bardsley2014,Wang2017}, but we generalizes it to Besov priors in a general parameter setting. The proposed RTO sampling method enables us to explore the sparsity-promoting properties of Besov priors in the samples of the posterior -- an aspect that, to the best of our knowledge, has not been previously studied in the literature.   

The main contributions of the paper are summarized as follows:
\begin{itemize}
    \item We introduce Besov priors with general parameters for linear Bayesian inverse problems.
    \item We develop an RTO-based method with a prior transform that can perform uncertainty quantification of a posterior arising from a linear Bayesian inverse problem with Besov priors under a general parameter setting.
    \item We show the validity and effectiveness of the proposed RTO method.
    \item We perform uncertainty quantification on multiple inverse problems with various choices of Besov priors, through numerical experiments. 
\end{itemize}

The remainder of the paper is organized as follows. In Section~\ref{sec2}, we
introduce Besov priors and explain how they are incorporated into linear
Bayesian inverse problems. In Section~\ref{sec3}, we present the RTO method,
construct a prior transformation map and combine it into the proposed RTO
method. In Section~\ref{sec4}, we do numerical experiments that showcase the
properties of Besov priors and the effectiveness of the proposed RTO method.
Finally, Section~\ref{sec5} concludes the paper with a summary and discussion.

\section{Linear Bayesian Inverse Problem with Besov Prior}
\label{sec2}

Let $X$ and $Y$ be separable Banach spaces, and let $A:X \to Y$ be a closed linear operator. Consider an ill-posed inverse problem, where an unknown quantity $f \in X$ needs to be recovered from a noisy measurement $y \in Y$ given by the equation
\begin{equation}
    y=A f +\epsilon, \label{eq:InvProb}
\end{equation}
where the noise $\epsilon:\Omega\rightarrow Y$ is a random variable understood in a measure-theoretic sense, see e.g., \citep{Folland1984}, with an underlying probability space $\left(\Omega, \Sigma, \mathbb{P}\right)$. We adopt the Bayesian approach to solve the inverse problem in \eqref{eq:InvProb} which allows for both regularization by imposing prior knowledge and quantification of uncertainties associated with measurement and modelling error. In the Bayesian approach, the unknown and the measurement are modelled as random variables, that is, $f:\Omega\rightarrow X$ and $y:\Omega\rightarrow Y$. We assign probability measures to $f$ and $y-Af$ with corresponding prior distribution $\pi_{\text{pr}}$  and likelihood distribution $\pi_{\text{like}}$, respectively. The solution to the inverse problem is then the conditional distribution of $f$ given data $y$, which is called the posterior distribution, and is denoted by $\pi_{\text{post}}$.

\subsection{Besov space random variables}
We motivate Besov random variables by briefly introducing Besov spaces and their properties. Besov spaces generalize well-known smoothness spaces such as H\"older and Sobolev spaces. They are a family of Banach spaces, where the regularity of functions can be measured accurately. Indeed, the regularity of functions in Besov spaces is measured through the function's partial derivatives and modulus of continuity \citep{cohen2003}. We define the $n$th order modulus of continuity for $f\in L^{p}(\mathbb{R}^{d})$ as
\begin{equation*}
    \omega_{n,p}(f,u)=\sup_{\substack{\|h\|_2\leq u}}\|\Delta_{h}^{n}f\|_{L^{p}(\mathbb{R}^{d})},
\end{equation*}
where $\Delta_{h}^{n}f(x)=\Delta_{h}^{n-1}(f(x)-f(x-h))$ is the $n$th order finite difference operator with the step-size $h$. A Besov space $B_{p,q}^{s}(\mathbb{R}^{d})$, with parameters $1\leq p,q<\infty$ and $s>0$, consists of functions $f\in L^{p}(\mathbb {R}^{d})$ with $\lfloor s \rfloor$ weak partial derivatives in $L^{p}(\mathbb{R}^{d})$ and where its modulus of continuity satisfies
\begin{equation}
    \{2^{sj}\omega_{n,p}(f,2^{-j})\}_{j\geq 0}\in \ell^{q}(\mathbb{N}),  \label{eq:BesovqCond}
\end{equation}
 for any $n>s$. The Besov space can be viewed as a Sobolev space $W^{s,p}(\mathbb{R}^{d})$ with the extra regularity condition, depending on $q$, given in \eqref{eq:BesovqCond} which fine-tunes the smoothness. 

The smoothness of functions in Besov spaces can be characterized by decay conditions on coefficients from expansions in a wavelet basis. Suppose $\phi\in L^{2}(\mathbb{R}^{d})$ generates a multiresolution analysis $(V_j)$ with dilation matrix $2I_d$. Let $\psi^{\ell}\in L^{2}(\mathbb{R}^{d})$, with $\ell\in \{1,2,\ldots,2^{d}-1\}$, denote the associated family of wavelets \citep{meyer_1993}. The function $\phi \in V_0$ is called the scaling function. We say that a  function $f$ is $r$-regular if $f\in C^{r}$ and
\begin{equation*}
    \left|\partial^{a}f(x) \right|\leq C_{l}\left(1+\|x\|_{1}\right)^{-l},
\end{equation*}
for any $l\in \mathbb{N}$ and any multi-index $\left|a\right|=a_{1}+a_{2}+\cdots + a_{d}\leq r$, where $C_{l}$ is a constant depending on $l$. A multiresolution analysis is $r$-regular if both the scaling function $\phi$ and the family of wavelets $\psi^{\ell}$ are $r$-regular. 
The dilation and translation of the wavelets and the scaling function are denoted as: 
\begin{align*}
    \psi^{\ell}_{j,k}(x)&=2^{jd/2}\psi^{\ell}(2^{j}x-k),\\
    \phi_{k}(x)&=\phi(x-k),
\end{align*}
with $\ell\in \{1,2,\ldots, 2^{d}-1\},\;j\geq 0,\; k\in \mathbb{Z}^{d}$.
Then any function $f\in L^{2}(\mathbb{R}^{d})$ has a wavelet expansion given by
\begin{equation}
        f=\sum_{k\in \mathbb{Z}^{d}}v_{k}\phi_{k}+\sum_{\ell=1}^{2^{d}-1}\sum_{j=0}^{\infty}\sum_{k\in \mathbb{Z}^{d}} w_{j,k}^{\ell}\psi_{j,k}^{\ell}, \label{eq:WaveletExpansion}
\end{equation}
where $v_{k}=\langle f,\phi_{k} \rangle_{L^{2}(\mathbb{R}^{d})}$ and $w_{j,k}^{\ell}$$=\langle f,\psi^{\ell}_{j,k} \rangle_{L^{2}(\mathbb{R}^{d})}$ with unconditional convergence in the norm. Under suitable regularity assumptions on the wavelet basis, membership in a Besov space can be characterized entirely by the decay of the wavelet coefficients in \eqref{eq:WaveletExpansion}. This wavelet characterization of Besov spaces is stated in the following proposition.

\begin{proposition}[Chapter 6, section 10 in \citep{meyer_1993}] \label{prop:WaveletBesChar}
Let $\psi^{\ell}, \ell\in \{1,2,\ldots,2^{d}-1\}$ be the family of wavelets and $\phi$ be the scaling function of a multiresolution analysis of $L^{2}(\mathbb{R}^{d})$ with a regularity $r\geq 1$. Let $1\leq p,q<\infty$ and $s>0$ such that $s<r$. Set $\kappa:=s+d/2-d/p$. The function $f\in L^{p}(\mathbb{R}^{d})$ has the wavelet expansion as in \eqref{eq:WaveletExpansion},
which satisfies 
\begin{equation*}
    \{v_{k}\}_{k\in \mathbb{Z}^{d}}\in \ell^{p}(\mathbb{Z}^{d}),
\end{equation*}
and
\begin{equation}
    \left\{2^{j\kappa}\left(\sum_{\ell=1}^{2^{d}-1}\sum_{k\in \mathbb{Z}^{d}}\left|w_{j,k}^{\ell}\right|^{p}\right)^{\frac{1}{p}}\right\}_{j\geq 0}\in \ell^{q}(\mathbb{N}),\label{eq:BesovCondition}
\end{equation}
if and only if $f\in B_{p,q}^{s}(\mathbb{R}^{d})$.
\end{proposition}
 
Proposition~\ref{prop:WaveletBesChar} characterizes Besov functions by their wavelet coefficients, and these conditions on the coefficients naturally
induce a norm on the Besov space $B_{p,q}^{s}(\mathbb{R}^{d})$ by
\begin{align*}
     \|f\|_{B_{p,q}^{s}(\mathbb{R}^{d})}=&\left(\sum_{k\in \mathbb{Z}^{d}}|v_{k}|^{p}\right)^{\frac{1}{p}}+\Biggl(\sum_{j=0}^{\infty}2^{j\kappa q} \nonumber\\ &\Biggl(\sum_{\ell=1}^{2^{d}-1}\sum_{k\in \mathbb{Z}^{d}}|w_{j,k}^{\ell}|^{p}\Biggr)^{\frac{q}{p}}\Biggr)^{\frac{1}{q}}, 
\end{align*}
where $\kappa=s+d/2-d/p$ \citep{triebel2006}.
 
The wavelet characterization of Besov functions in Proposition~\ref{prop:WaveletBesChar} is fundamental to the definition of Besov space random variables. Building on a similar characterization, \citep{Dashti_Stuart2012,Lassas2009} construct random variables on Besov spaces of 1-periodic functions defined over the torus $\mathbb{T}^{d}$, assuming for simplicity that the space parameters satisfy $p=q$. With these assumptions, the definition of Besov space priors is given in the following.

 \begin{definition}[\cite{Dashti_Stuart2012}] \label{def:BesovRand}
Let $\psi^{\ell}, \ell\in \{1,2,\ldots,2^{d}-1\}$ be the family of wavelets of a multiresolution analysis in $L^{2}(\mathbb{R}^{d})$ with regularity $r\geq 1$. Let $1\leq p<\infty$ and $s>0$ such that $s<r$. Let $\mathcal{K}_{j}=\{0,\ldots,2^{j}-1\}^{d}$ be an index set and $\{\xi_{j,k}^{\ell}\}_{j\geq 0,k\in \mathcal{K}_{j}}, \xi_{0}$ be i.i.d. real valued, generalized Gaussian  random variables with density $\pi_{\xi}(x)\propto \exp (-\tau|x|^{p})$ with $\tau>0$. Let $f$ be defined as
\begin{align}
    f(x)=\xi_{0}+\sum_{\ell=1}^{2^{d}-1}\sum_{j=0}^{\infty}\sum_{k\in \mathcal{K}_{j}}&2^{-j\kappa}\xi^{\ell}_{j,k}\psi^{\ell}_{j}(x-k2^{-j}), \label{eq:BesovExpansion}
\end{align}
for a.e. $x\in \mathbb{T}^{d}$, where 
\[ 
\psi_{j}^{\ell}(x)=2^{jd/2}\sum_{m\in \mathbb{Z}^{d}}\psi^{\ell}(2^{j}(x-m))
\] 
is the 1-periodization of $\psi_{j,k}^{\ell}$. We say that $f$ is a Besov $B_{p,p}^{s}(\mathbb{T}^{d})$ random variable if the series in \eqref{eq:BesovExpansion} converges.
\end{definition}
The norm on the periodic Besov space $B_{p,p}^{s}(\mathbb{T}^{d})$ is the $\ell^{p}(\mathbb{N})$ norm on the sequence in \eqref{eq:BesovCondition} with the periodic restriction, that is
\begin{align*}
    \|f\|_{B_{p,p}^{s}(\mathbb{T}^{d})}=\biggl(&\left|v_{0}\right|^{p}+\sum_{\ell=1}^{2^{d}-1}\sum_{j=0}^{\infty}\sum_{k\in \mathcal{K}_{j}}\nonumber \\ &2^{j\kappa p}|w^{\ell}_{j,k}|^{p}\biggr)^{\frac{1}{p}},
\end{align*}
where $w^{\ell}_{j,k}=\langle f,\psi^{\ell}_{j}(x-k2^{-j}) \rangle_{L^{2}(\mathbb{T}^{d})}$.
The probability density function for a Besov space random variable following Definition~\ref{def:BesovRand} is formally given by $\pi_{\text{Besov}}(f)\propto \exp(-\tau\|f\|^{p}_{B_{p,p}^{s}(\mathbb{T}^{d})})$. 
\begin{proposition}[\cite{Lassas2009}] The Besov random variable $f$ in Definition~\ref{def:BesovRand} takes values in the Besov space $B_{p,p}^{t}(\mathbb{T}^{d})$ if and only if the parameters satisfy $t<s-\frac{d}{p}$.
\end{proposition}

\begin{remark}
The definition of Besov spaces $B_{p,q}^{s}(\mathbb{R}^{d})$ can be extended to $s\in \mathbb{R}$ where the elements are now in the space of tempered distributions $\mathcal{S}'\left(\mathbb{R}^{d}\right)$. The conditions for membership of a Besov space $B_{p,q}^{s}(\mathbb{R}^{d})$ is, for example, given in \citep{Sawano2018} as Definition 2.1. Additionally, with this new definition of Besov spaces, Proposition~\ref{prop:WaveletBesChar} and Definition~\ref{def:BesovRand} can also be extended to $s\in \mathbb{R}$ which is done in \citep{meyer_1993} and \citep{Lassas2009}, respectively.
\end{remark}

\subsection{Discretization of the Bayesian inverse problem}

In this subsection, we discretize the prior and the likelihood distributions to obtain the discrete posterior. We discretize the inverse problem to a finite form given by
\begin{equation}
    y_{m}=Af_{n}+\epsilon_{m}, \label{eq:DiscreteInverseProblem}
\end{equation}
where the measurement $y_{m}$, the unknown $f_{n}$, and the noise $\epsilon_{m}$ are random variables that take values in $\mathbb{R}^{m}$, $\mathbb{R}^{n}$, and $\mathbb{R}^{m}$, respectively. The linear forward map $A$ is now represented by a matrix $A\in\mathbb{R}^{m\times n}$. We assume that the noise $\epsilon_{m}$ is Gaussian and distributed according to $\epsilon_{m}\sim N(0,\sigma^{2}I_{m})$ with $\sigma>0$. Thus, the likelihood is Gaussian and takes the form
\begin{equation}
    \pi_{\text{like}}(y_{m}|f_{n})\propto \exp\left(-\frac{1}{2\sigma^{2}}\|Af_{n}-y_{m}\|_{2}^{2}\right). \label{eq:DiscreteLikelihood}
\end{equation}
The discrete Besov random variable is obtained by truncating the random expansion in \eqref{eq:BesovExpansion} at some finite integer $J_{\text{max}}\geq1$, that is
\begin{align}
    f_{\text{dis}}(x)=\xi_{0}+\sum_{\ell=1}^{2^{d}-1}&\sum_{j=0}^{J_{\text{max}}-1}\sum_{k\in \mathcal{K}_{j}}2^{-j\kappa}\xi^{\ell}_{j,k}\nonumber \\ &\psi^{\ell}_{j}(x-k2^{-j}),\label{eq:FiniteBesovExpansion}
\end{align}
for a.e. $x\in \mathbb{T}^{d}$. 
The norm of the truncated Besov random variable is
\begin{align}
    \label{eq:DiscreteBesovNorm}
 \|f_{\text{dis}}\|_{B_{p,p}^{s}(\mathbb{T}^{d})}=\biggl(&\left|v_{0}\right|^{p}+\nonumber \sum_{\ell=1}^{2^{d}-1}\sum_{j=0}^{J_{\text{max}}-1}\sum_{k\in \mathcal{K}_{j}}\\& 2^{j\kappa p}|w^{\ell}_{j,k}|^{p}\biggr)^{\frac{1}{p}}.
\end{align}

We model the unknown $f_{n}$ as a discrete Besov random function $f_{\text{dis}}$ on an equidistant grid in $[0,1[^{d}$ with interspacing $2^{-J_{\text{max}}}$. Let $\delta_{n}$ be an approximation of the wavelet coefficients $w_{j,k}^{\ell}$ associated with $f_{\text{dis}}$. If $n=2^{J_{\text{max}}d}$ then $\delta_{n}$ can be computed from the unknown $f_{n}$ using 
the periodic fast wavelet transform algorithm \citep{Mallat1989}. We denote the matrix representation of this linear transformation by $W\in\mathbb{R}^{n\times n}$, i.e., 
\begin{equation*}
    \delta_{n}=Wf_{n}.
\end{equation*}
We define a diagonal matrix $S\in \mathbb{R}^{n\times n}$ of Besov weights from \eqref{eq:DiscreteBesovNorm} as 
\begin{align*}
    S_{1,1}=1&,\quad S_{i,i}=2^{j\kappa},
\end{align*}
for $2^{jd}+1\leq i \leq 2^{(j+1)d}$ and for $j=0,1,\ldots, J_{\text{max}}-1$. 
 Now, the discrete Besov norm can be approximated in terms of $S$ and $W$ as
\begin{equation*}
   \|f\|_{B_{p,p}^{s}(\mathbb{T}^{d})}\approx \|SWf_{n}\|_{p}.
\end{equation*}
Setting $B=SW\in\mathbb{R}^{n\times n}$, the discrete Besov prior can be written as
\begin{equation}
    \pi_{\text{pr}}(f_{n})\propto \exp\left(-\tau\|Bf_{n}\|_{p}^{p}\right). \label{eq:DiscreteBesovPrior}
\end{equation}
The discrete posterior of the Bayesian inverse problem in \eqref{eq:DiscreteInverseProblem} can now be computed by combining the discrete likelihood and the discrete prior as
\begin{align}
    \pi_{\text{post}}(f_{n}|y_{m})\propto \exp\biggl(&-\frac{1}{2\sigma^{2}}\|Af_{n}-y_{m}\|_{2}^{2}\nonumber\\&-\tau\|Bf_{n}\|_{p}^{p}\biggr). \label{eq:DiscretePosterior}
\end{align}

\section{Randomize-Then-Optimize sampling method}
\label{sec3}

In this section, we introduce the RTO method for exploring the posterior distribution in \eqref{eq:DiscretePosterior}. A common approach to explore such posteriors is using MCMC methods which approximates the posterior with a set of samples. These methods comprise a proposal step of proposing correlated samples based on previously obtained samples, followed by an acceptance/rejection step which compares and selects samples, based on their relative posterior probability. However, such methods are inefficient, i.e., most samples are rejected, when the prior distribution of the posterior is a generalized Gaussian, particularly with $p<2$ \citep{chen2018robust,lucka2016fast,Wang2017}.

RTO is a sampling method that draws samples from a posterior distribution where the negative-log-posterior has a least squares form with Tikhonov regularization. This method approximates the posterior with a set of samples drawn by solving a randomly perturbed least squares problem with Tikhonov regularization. Least squares problems with Tikhonov regularization have been extensively studied in the past decades and efficient numerical methods are available.

In the linear case, the posterior is Gaussian, and the RTO method draws exact independent samples from the posterior. Unlike many MCMC methods, the proposed sample in the RTO method is independent of the previous samples, yielding a higher probability of acceptance and lower sample autocorrelation. This makes RTO an effective proposal-based MCMC method that we want to exploit in the setting of the posterior in \eqref{eq:DiscretePosterior}.

We start by introducing the RTO method. The RTO method is not directly applicable to the posterior in \eqref{eq:DiscretePosterior} so we introduce a transformation that enables us to change variables in the posterior such that the transformed posterior becomes applicable for RTO.

\subsection{RTO for a linear inverse problem with Gaussian prior}
The RTO method draws samples from a posterior of the form
\begin{equation}
    \pi_{\text{post}}(f_{n}|y_{m}) \propto \exp\left(-\frac{1}{2}\|F(f_{n})\|_{2}^{2}\right),\label{eq:GaussPost}
\end{equation}
where the map $F:\mathbb{R}^{n}\rightarrow \mathbb{R}^{m+n}$ combines a Gaussian likelihood and a Gaussian prior. In the case of an inverse problem of the form \eqref{eq:DiscreteInverseProblem} with a Gaussian likelihood as in \eqref{eq:DiscreteLikelihood}, the mapping $F$ takes the form
\begin{equation*}
    F(f_{n})=\begin{bmatrix}L^{-1/2}f_{n}\\ \sigma^{-1}\left(Af_{n}-y_{m}\right) \end{bmatrix},
\end{equation*}
where $L$ is positive semidefinite and $f_{n}\sim N(0,L)$. The RTO method suggests that the samples from this distribution can be achieved by solving the following stochastic least squares problem (see Prop.~\ref{prop:RTO})
\begin{equation}
    f_{\text{sample}}=\argmin_{\substack{f_{n}}}\left\|F(f_{n})-v_{n+m}\right\|^{2}_{2},\label{eq:StochasticLeastSquares}
\end{equation}
where $v_{n+m}\sim N(0,I_{m+n}).$
We can achieve computationally stable solutions to the stochastic optimization problem in \eqref{eq:StochasticLeastSquares} using the thin QR decomposition of the system matrix 
\begin{equation*}
\mathcal{A}=\begin{bmatrix}L^{-1/2}\\ \sigma^{-1}A\end{bmatrix}=QR,
\end{equation*}
where $Q\in \mathbb{R}^{(m+n)\times n}$ has orthonormal columns, and $R\in \mathbb{R}^{n\times n}$ is upper triangular. Suppose $\mathcal{A}$ has full rank. Then $R$ is invertible, and the solution to \eqref{eq:StochasticLeastSquares} is unique and can be expressed in terms of a stochastic equation
\begin{equation}
    Q^{T}\mathcal{A}f_{\text{sample}}=Q^{T}\left(\begin{bmatrix}0_n\\ \sigma^{-1}y_{m}\end{bmatrix}+v_{m+n}\right).\label{eq:RTOequation}
\end{equation}
%where $v_{m+n}\sim N(0,I_{m+n})$.
We can use the stochastic equation \eqref{eq:RTOequation}, together with the following proposition, to determine which distribution the RTO samples come from.

\begin{proposition}[\cite{bardsley2018computational}]\label{prop:RTO}
Let $G:\mathbb{R}^{N}\rightarrow \mathbb{R}^{N}$ be a one-to-one function with continuous first partial derivatives and invertible Jacobian matrix $J(x)$ for all $x\in \mathbb{R}^{N}$. Then the random vector $x$ is defined through
\begin{equation*}
    G(x)=v,\quad v\sim N(\mu,\Sigma),
\end{equation*}
has the probability density function
\begin{equation*}
    \pi(x)\propto \left|J(x)\right|\exp\left(-\frac{1}{2}\|G(x)-\mu\|_{2,\Sigma^{-1}}^{2}\right),
\end{equation*}
where $|\cdot|$ denotes the absolute value of the determinant and $\|x\|_{2,\Sigma^{-1}}^{2}=x^{T}\Sigma^{-1}x$.
\end{proposition}
 
To see that $f_{\text{sample}}$ comes from the posterior distribution in \eqref{eq:GaussPost}, we define the map $G=Q^{T}\mathcal{A}$ and use Proposition~\ref{prop:RTO} with $v\sim N(Q^{T}\left[0_n,\sigma^{-1}y_{m}\right]^{T},I_{n})$ to define the RTO distribution $\pi_{\text{RTO}}$.  The resulting RTO distribution can easily be identified as the Gaussian distribution $\pi_{\text{RTO}}(f_{n})\propto \exp\left(-\frac{1}{2}\|Q^{T}F(f_{n})\|_{2}^{2}\right)$. Additionally, it can be shown using $QQ^{T}\mathcal{A}=\mathcal{A}$ that $\exp\left(-\frac{1}{2}\|Q^{T}F(f_{n})\|_{2}^{2}\right)\propto \exp\left(-\frac{1}{2}\|F(f_{n})\|_{2}^{2}\right)$. We conclude that  the RTO samples are exact in the linear case.

\subsection{Prior transformation}

Our goal is to sample from the posterior in \eqref{eq:DiscretePosterior} using the RTO method, but for $p\neq 2$ the posterior is not on the form as in \eqref{eq:GaussPost} so RTO does not apply. To overcome this problem, we introduce a prior transform from Besov to standard Gaussian using an inverse cumulative distribution function (inverse CDF) type method. A generalised Gaussian distribution extends the classical Gaussian and Laplace families \citep{Nadarajah2005}. Its probability density function (PDF) is given by 
\begin{equation*}
    \pi_{\alpha,p}(x)=\frac{p}{\alpha \Gamma(\frac{1}{p})}\exp\left(-\left|\frac{x-\mu}{\alpha}\right|^{p}\right),
\end{equation*}
where $\alpha$ is the scale, $p$ is the shape, $\mu$ is the mean, and $\Gamma(\cdot)$ is the gamma function.
We choose the parameters as $\mu=0$ and $\alpha=\left(\Gamma(1/p)/\Gamma(3/p)\right)^{1/2}\lambda^{-1/p}$ with $\lambda>0$ that relates to the prior parameter $\tau$ in \eqref{eq:DiscreteBesovPrior} by $\left(\Gamma(1/p)/\Gamma(3/p)\right)^{-p/2}\lambda=\tau$. With these choices, $\lambda$ is a scaled regularization parameter that controls the variance such that when $\lambda=1$ the variance is 1.
We introduce an invertible mapping $g_{\text{1D}}: \mathbb{R} \to \mathbb{R}$ from the standard Gaussian random variable to a generalized Gaussian random variable using an inverse CDF method. We set
\begin{equation*}
    g_{\text{1D}}(h)=\Phi_{\alpha,p}^{-1}(\phi_{G}(h)),
\end{equation*}
where $\Phi_{\alpha,p}$ is the CDF of a generalized Gaussian and $\phi_{G}$ is the CDF of a standard Gaussian.
Since $f_{n}$ is distributed according to a Besov prior, each element of $Bf_{n}$ is i.i.d. according to a generalized Gaussian. Consequently, we can construct an invertible mapping $g:\mathbb{R}^{n}\rightarrow \mathbb{R}^{n}$ that maps a standard Gaussian vector $h_{n}$ to $Bf_{n}$ using $g_{\text{1D}}$ component wise
\begin{align*}
    Bf_{n}&=g(h_{n})\nonumber\\&=\left[g_{\text{1D}}((h_{n})_{1}),\ldots,g_{\text{1D}}((h_{n})_{n}) \right]^{T}. 
\end{align*}
With this construction, we can introduce a prior transform from Besov prior $f_{n}$ to standard Gaussian prior $h_{n}$ by
\begin{equation}
T(h_{n})=B^{-1}g(h_{n})=f_{n}. \label{eq:PriorMap}
\end{equation}
Thus, we map the posterior in \eqref{eq:DiscretePosterior} into the following posterior with respect to the new variable $h_{n}$:
\begin{align}
    \pi_{\text{post}}(h_{n}|y_{m}) \propto \exp\biggl(&-\frac{1}{2\sigma^{2}}\|AB^{-1}g(h_{n})-y_{m}\|_{2}^{2}\nonumber\\&-\frac{1}{2}\|h_{n}\|_{2}^{2}\biggr). \label{eq:TransformedPosterior}
\end{align}
Notice that the forward operator $AB^{-1}g$ that maps $h_{n}\mapsto y_{m}$ is now nonlinear since $g$ is nonlinear.

\subsection{RTO sampling for Linear inverse problem with Besov prior}
The transformed posterior \eqref{eq:TransformedPosterior} is on the form \eqref{eq:GaussPost} with 
\begin{equation*}
    F(h_{n})=\begin{bmatrix} h_{n}\\ \sigma^{-1}\left(AB^{-1}g(h_{n})-y_{m}\right)\end{bmatrix},
\end{equation*}
and is applicable to the RTO method. In addition, we can transform a sample from the posterior in \eqref{eq:TransformedPosterior} to the desired posterior \eqref{eq:DiscretePosterior} using the invertible prior transformation \eqref{eq:PriorMap}. The transformed posterior \eqref{eq:TransformedPosterior} gives rise to a stochastic equation on the form \eqref{eq:RTOequation} where $\mathcal{A}(h_{n})=\left[h_{n}^T,\sigma^{-1}(AB^{-1}g(h_{n}))^T\right]^{T}$ is nonlinear. In this scenario, an invertible $Q$ can be chosen freely, but such that the mapping $G(\cdot)\coloneqq Q^{T}\mathcal{A}(\cdot)$ satisfies the conditions of Proposition~\ref{prop:RTO} which results in the following RTO distribution 
\begin{equation*}
    \pi_{\text{RTO}}(h_{n})\propto \left|Q^{T}J_{\mathcal{A}}(h_{n}) \right|\exp\left(-\frac{1}{2}\|Q^{T}F(h_{n})\|_{2}^{2}\right).
\end{equation*}
 We can rewrite the RTO distribution $\pi_{\text{RTO}}$ in terms of the target posterior in \eqref{eq:TransformedPosterior} as
\begin{equation*}
    \pi_{\text{RTO}}(h_{n})\propto c(h_{n})\exp\left(-\frac{1}{2}\|F(h_{n})\|_{2}^{2}\right),
\end{equation*}
where
\begin{align}
    c(h_{n})=\left|Q^{T}J_{\mathcal{A}}(h_{n})\right|\exp\biggl(&\frac{1}{2}\|F(h_{n})\|_{2}^{2}\nonumber\\-&\frac{1}{2}\|Q^{T}F(h_{n})\|_{2}^{2}\biggr). \label{eq:MHweights}
\end{align}

Let $h_{\text{MAP}}$ be the maximum a posteriori (MAP) estimate of the posterior \eqref{eq:TransformedPosterior}, that is,
\begin{equation}
    h_{\text{MAP}}=\argmin_{\substack{h_{n}}}\|F(h_{n})\|_{2}^{2}.\label{eq:MAPRTO}
\end{equation}
In \citep{Bardsley2014,bardsley2018computational} they suggest picking $Q$ from the thin QR-factorization of the Jacobian of $\mathcal{A}$ evaluated at $h_\text{MAP}$, and we will likewise use this choice. In other words, if $J_{\mathcal{A}}$ is the Jacobian matrix of $\mathcal{A}$ then we compute Q from the QR-factorization 
\begin{equation*}
    J_{\mathcal{A}}(h_{\text{MAP}})=\begin{bmatrix}Q & \Tilde{Q} \end{bmatrix}\begin{bmatrix}R \\0\end{bmatrix}.
\end{equation*}
With this choice, $QQ^{T}$ is an orthogonal projection guiding the samples of $\pi_{\text{RTO}}$ towards regions of higher probability of the target posterior \eqref{eq:TransformedPosterior}.

To compute the RTO samples of $\pi_{\text{RTO}}$ we solve the nonlinear stochastic optimization problem given by
\begin{equation}
    h_{\text{prop}}=\argmin_{\substack{h_{n}}}\|Q^{T}\left(F(h_{n})-v_{m+n}\right)\|_{2}^{2}.\label{eq:NonlinearStochasticLSQ}
\end{equation}
%where $v_{m+n}\sim N(0,I_{m+n}).$
To ensure that $h_{\text{prop}}$ satisfies the stochastic equation \eqref{eq:RTOequation}, which makes it a sample from $\pi_{\text{RTO}}$, we need to check that the cost function of the optimization problem is sufficiently small, i.e.,
\begin{equation}
    \|Q^{T}\left(F(h_{\text{prop}})-v_{m+n}\right)\|_{2}^{2}<\eta, \label{eq:Costinequality}
\end{equation}
with $\eta$ small.
Since $\pi_{\text{RTO}}\neq \pi_{\text{post}}$, the samples must be corrected using the Metropolis-Hastings (MH) acceptance-rejection technique, and we use $h_{\text{prop}}$ as a proposal. This ensures that the accepted RTO samples belong to the desired posterior in \eqref{eq:TransformedPosterior}.

To illustrate the MH acceptance-rejection step, suppose that $h^{i-1}_{\text{sample}}$ is the $(i-1)$th sample in the RTO method. We compute $h^{i}_{\text{prop}}$ from \eqref{eq:NonlinearStochasticLSQ} and then compute the acceptance probability
\begin{align*}
    p_{\text{MH}}^{i}&=\frac{\pi_{\text{post}}(h_{\text{prop}}^{i}|y_{m}) \pi_{\text{RTO}}(h_{\text{sample}}^{i-1})}{\pi_{\text{post}}(h_{\text{sample}}^{i-1}|y_{m}) \pi_{\text{RTO}}(h_{\text{prop}}^{i})}\nonumber\\&=\frac{c(h_{\text{sample}}^{i-1})}{c(h_{\text{prop}}^{i})}.
\end{align*}
Note that to evaluate $\pi_{\text{RTO}}$ we must evaluate the Jacobian with respect to $h_{n}$
\begin{equation}
    J_{\mathcal{A}}(h_{n})=\begin{bmatrix}I_{n}\\ \sigma^{-1}AB^{-1}J_{g}(h_{n})\end{bmatrix}. \label{eq:TransformJacobian}
\end{equation}
The Jacobian of the mapping is diagonal and given by
\begin{equation}
    J_{g}(h_{n})=\begin{bmatrix}g_{\text{1D}}^{'}((h_{n})_{1}) & & & \\
    %& g_{\text{1D}}^{'}((h_{n})_{2}) & &  \\ 
    & & \ddots & \\
    & & & g_{\text{1D}}^{'}((h_{n})_{n}) \end{bmatrix},\label{eq:Jacobiangmap}
\end{equation}
and 
\begin{equation*}
    g^{'}_{\text{1D}}(h)=\frac{\pi_{\text{G}}(h)}{\pi_{\alpha,p}[\Phi_{\alpha,p}^{-1}(\phi_{G}(h))]},
\end{equation*}
where $\pi_{\text{G}}$ is the PDF of the standard Gaussian distribution.

Finally, we summarize our discussion in Alg.~\ref{alg:RTOMHBesov}, which describes the RTO-MH algorithm used to sample from the posterior \eqref{eq:DiscretePosterior}.

\begin{algorithm}[t] 
\caption{RTO-MH sampling from a posterior with a Besov prior}\label{alg:RTOMHBesov}
\begin{algorithmic}[1]
\REQUIRE Initial guess $h_{\text{sample}}^{0}$, $\eta>0$ in \eqref{eq:Costinequality}, and the amount of samples $n_{\text{samp}}$.\\
\COMMENT{Compute line 1:4 offline;}
\STATE Specify the prior transform \eqref{eq:PriorMap} together with the Jacobian \eqref{eq:Jacobiangmap}.
\STATE Compute the MAP estimate $h_{\text{MAP}}$ of the transformed posterior in \eqref{eq:MAPRTO}
\STATE Compute $Q$ from the thin QR factorization of $J_{\mathcal{A}}(h_{\text{MAP}})$ where $J_{\mathcal{A}}$ is given in \eqref{eq:TransformJacobian}.\\
\STATE Compute $c(h_{\text{sample}}^{0})$ as in \eqref{eq:MHweights} \\
\COMMENT{RTO procedure lines 5:9 computed online;}
\FOR{$i=1,\ldots,n_{\text{samp}}$}
\STATE Draw a sample $v_{n+m}^{i}$ from $N(0,I_{n+m})$
\STATE Compute the proposal sample $h_{\text{prop}}^{i}$ by solving \eqref{eq:NonlinearStochasticLSQ} and compute the cost function by $e^{i}=\left\|Q^{T}\left(F(h_{\text{prop}}^{i})-v_{n+m}^{i}\right)\right\|_{2}^{2}$. 
\STATE Compute the MH coefficients $c(h_{\text{prop}}^{i})$ according to eq.~\eqref{eq:MHweights}.
\ENDFOR\\
\COMMENT{MH acceptance-rejection procedure lines 10:17 computed online};
\FOR{$i=1,\ldots, n_{\text{samp}}$}
\STATE Draw a sample $u$ from a uniform distribution on $[0,1]$.
\IF{$u<c(h_{\text{sample}}^{i-1})/c(h_{\text{prop}}^{i})$ and $e^{i}<\eta$}
\STATE $h_{\text{sample}}^{i}=h_{\text{prop}}^{i}$.
\ELSE
\STATE $h_{\text{sample}}^{i}=h_{\text{sample}}^{i-1}$.
\ENDIF
\ENDFOR\\
\COMMENT{Back transformation lines 18:20 computed online};
\FOR{$i=1,\ldots, n_{\text{samp}}$}
\STATE Transform the samples back to the non-Gaussian case by computing
\[
f_{\text{sample}}^{i}=T(h_{\text{sample}}^{i}).
\]
\ENDFOR
\end{algorithmic}
\end{algorithm}

\section{Numerical Tests}
\label{sec4}

In this section, we apply the proposed RTO-MH algorithm to three inverse problems with various Besov priors namely the 1D inpainting, the 1D deconvolution, and the 2D sparse angle computed tomography. The numerical tests are used to illustrate the properties of the Besov prior, to evaluate the effectiveness of the RTO-MH method, and to investigate the capabilities of the RTO-MH method. All the code used for producing the numerical tests, results, and visualizations are publicly available in the following GitHub repository \url{https://github.com/AndreasHorst/BP_RTO}.

\subsection{1D Inpainting} 
\label{sec:Inpainting}

Inpainting is the task of recovering a function $f$ on a domain $\Omega$ where it is unknown in a subdomain $D\subset \Omega$ called the inpainting region. We consider inpainting on the torus $\Omega=\mathbb{T}$ with inpainting region $D\subset \mathbb{T}$. The inpainting problem can be modeled as an inverse problem with forward operator $f \mapsto f\vert_{\mathbb{T}\setminus D}$, where $f\vert_{\mathbb{T}\setminus D}$ is the restriction of $f$ to the complement of $D$ in $\mathbb{T}$. 
% \mathcal{X}_{\mathbb{T}\setminus D}$, where
% \begin{equation}
%     \mathcal{X}_{\mathbb{T}\setminus D}(x)=\begin{cases}&1 \quad x\in \mathbb{T}\setminus D, \\ & 0 \quad x\in D. \end{cases}
% \end{equation}

\begin{figure}[h]
     \centering
     \includegraphics[width=76.0mm]{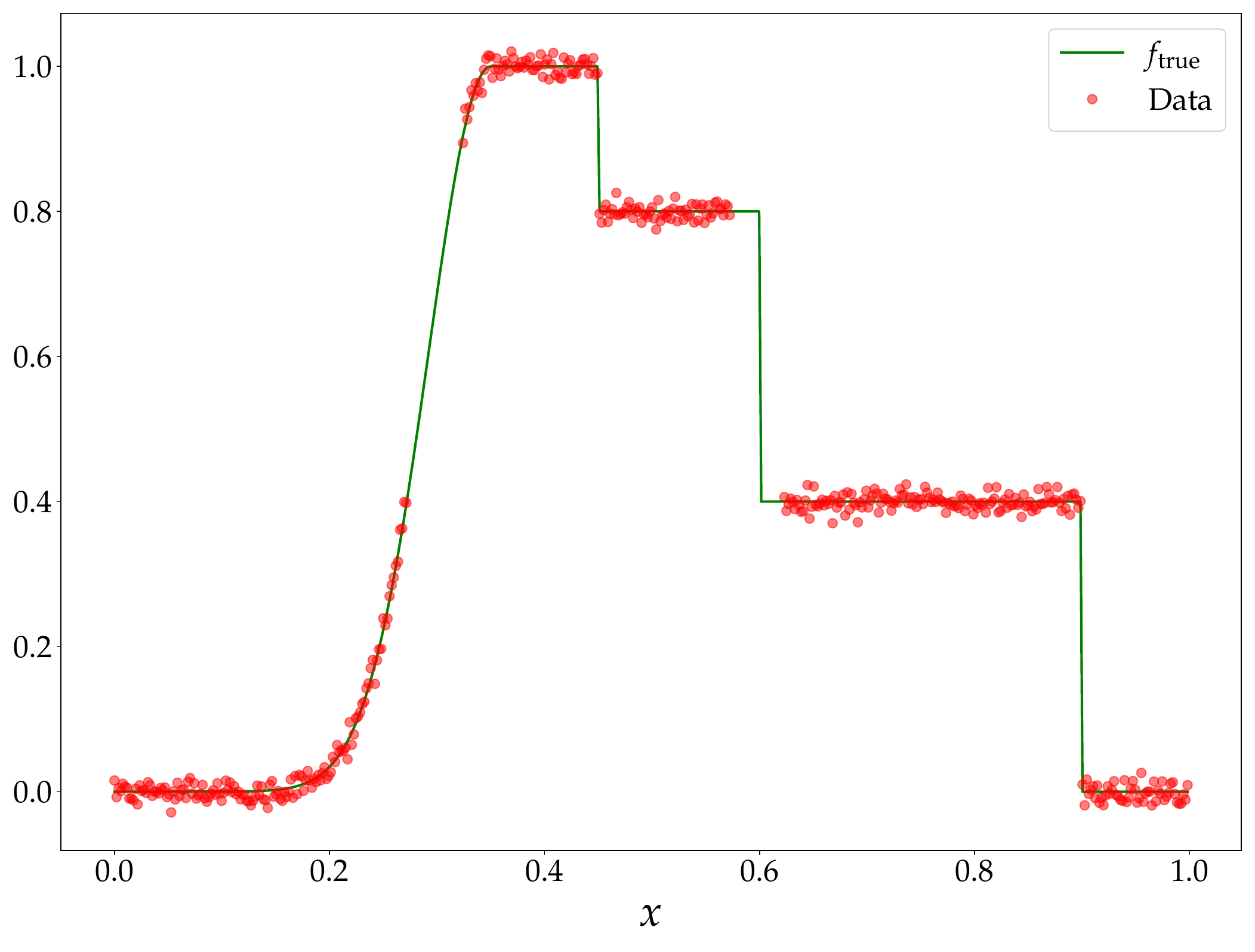}
     \caption{Ground truth signal and noisy data for the 1D inpainting problem.}
     \label{fig:InpaintingData}
 \end{figure}

 In the discrete case in 1D, we formulate the inpainting problem as
 \begin{equation*}
     y_{m}= \tilde{I}_{m}f_{n}+\epsilon_{m},
 \end{equation*}
 where $m<n$ and $\tilde{I}_{m}\in \mathbb{R}^{m\times n}$ is the identity matrix $I_{n}\in \mathbb{R}^{n\times n}$ with $n-m$ rows removed.
 We use this 1D inpainting problem to illustrate the properties of the Besov prior. The ground truth is chosen to be a signal in $[0,1]$ which consists of piecewise constant and piecewise smoothly varying parts. We discretize it on the grid of points given by $x_{n}=\frac{i}{n}$, $i=0,1,\ldots, n-1$ with $n=512$. The forward operator removes the parts of the signal in the intervals $[0.1,0.15]$ and $[0.425,0.475]$, and it results in $m=461$ data elements in $y_{m}$. The noisy data $y_{m}$ is computed by adding centered i.i.d. Gaussian noise to the data, where the standard deviation $\sigma$ is chosen such that $\frac{\|\epsilon_{m}\|_{2}}{\|f_{\text{true}}\|_{2}}=0.02$. The ground truth and the noisy data are shown in Figure~\ref{fig:InpaintingData}.
 
 \subsubsection{Influence of the wavelet basis}
 \label{sec:ChangeWavelet}
\begin{figure}[t]
     \centering
     \includegraphics[width=76.0mm]{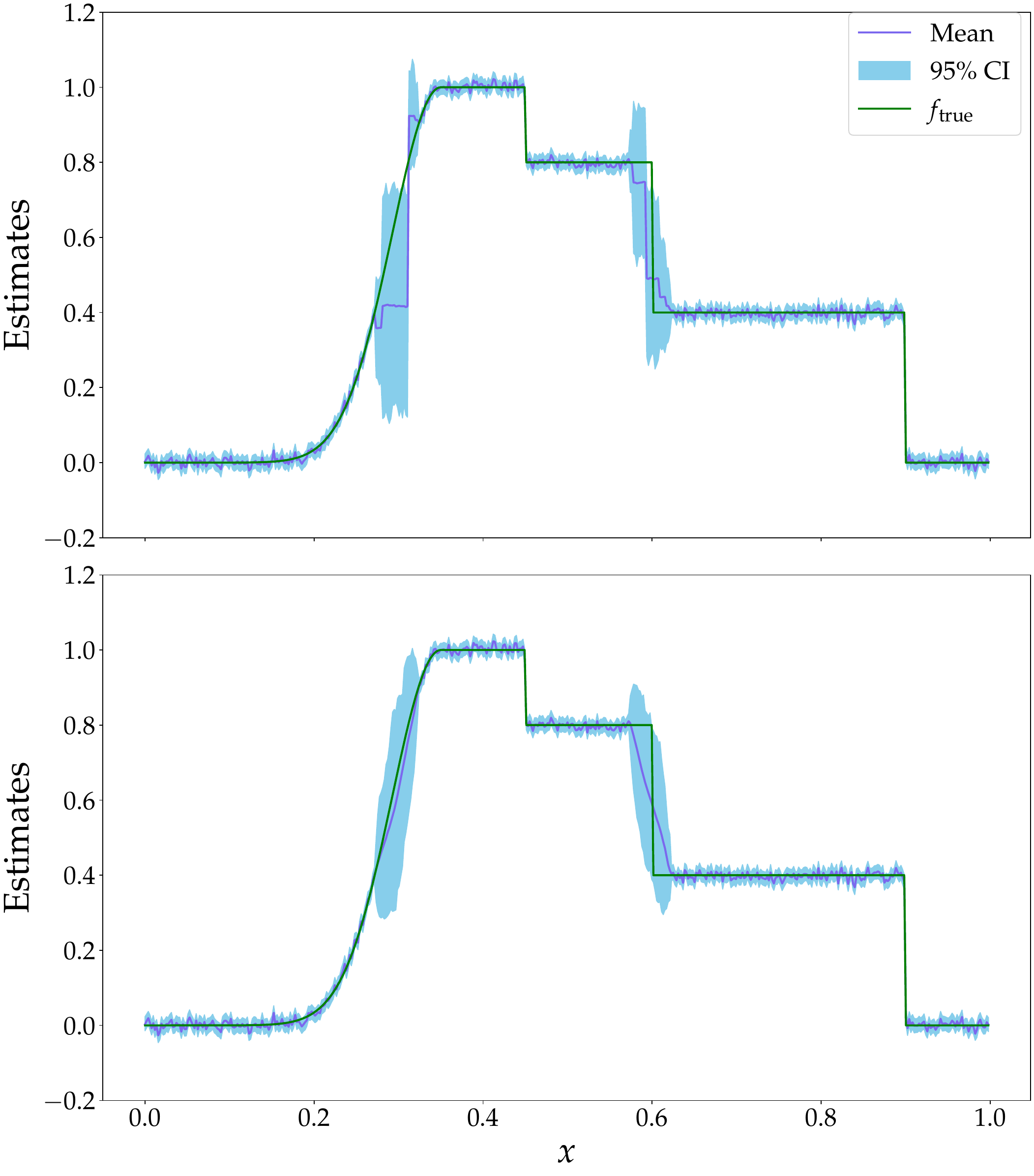}
     \caption{Posterior estimates using the Haar wavelet Besov prior (in the top) and the db8 wavelet Besov prior (in the bottom) with parameters $s=1.2$ and $p=1.5$.}
     \label{fig:Haar_db8_Inpaiting}
 \end{figure}

 The first inpainting test scenario investigates how the choice of wavelet basis in the Besov prior influences the posterior's ability to represent smooth and discontinuous functions. We perform the test with fixed parameters $s=1.2$, $p=1.5$, $\lambda=0.025$, and compare the results from using the Haar wavelet and the Daubechies-8 (db8) wavelet \citep{Daubechies1988Orthonormal}. In this test, we generate $10^4$ posterior samples from which we use the accepted samples to compute the posterior mean estimate and the 95\% credible interval (CI) estimate. The test resulted in 4723 and 4726 accepted samples in the Haar wavelet and db8 wavelet case, respectively. The resulting posterior estimates are displayed in Figure~\ref{fig:Haar_db8_Inpaiting}.
 
 In the inpainting regions, we can clearly observe the influence of the choice of wavelet basis in the Besov prior since the prior is the dominant term of the posterior in these regions. With the Haar wavelet, the Besov prior promotes discontinuous jumps in the posterior mean estimate, while with the db8 wavelet, the Besov prior promotes a smooth posterior mean estimate. Consequently, the smoothness of the wavelet basis in the Besov prior determines the smoothness of the posterior mean estimate. Additionally, we notice that the posterior 95\% CI estimate does not depend significantly on the choice of wavelet basis in the Besov prior.

 \subsubsection{Influence of the parameters}
\begin{figure}[t]
     \centering
     \includegraphics[width=76.0mm]{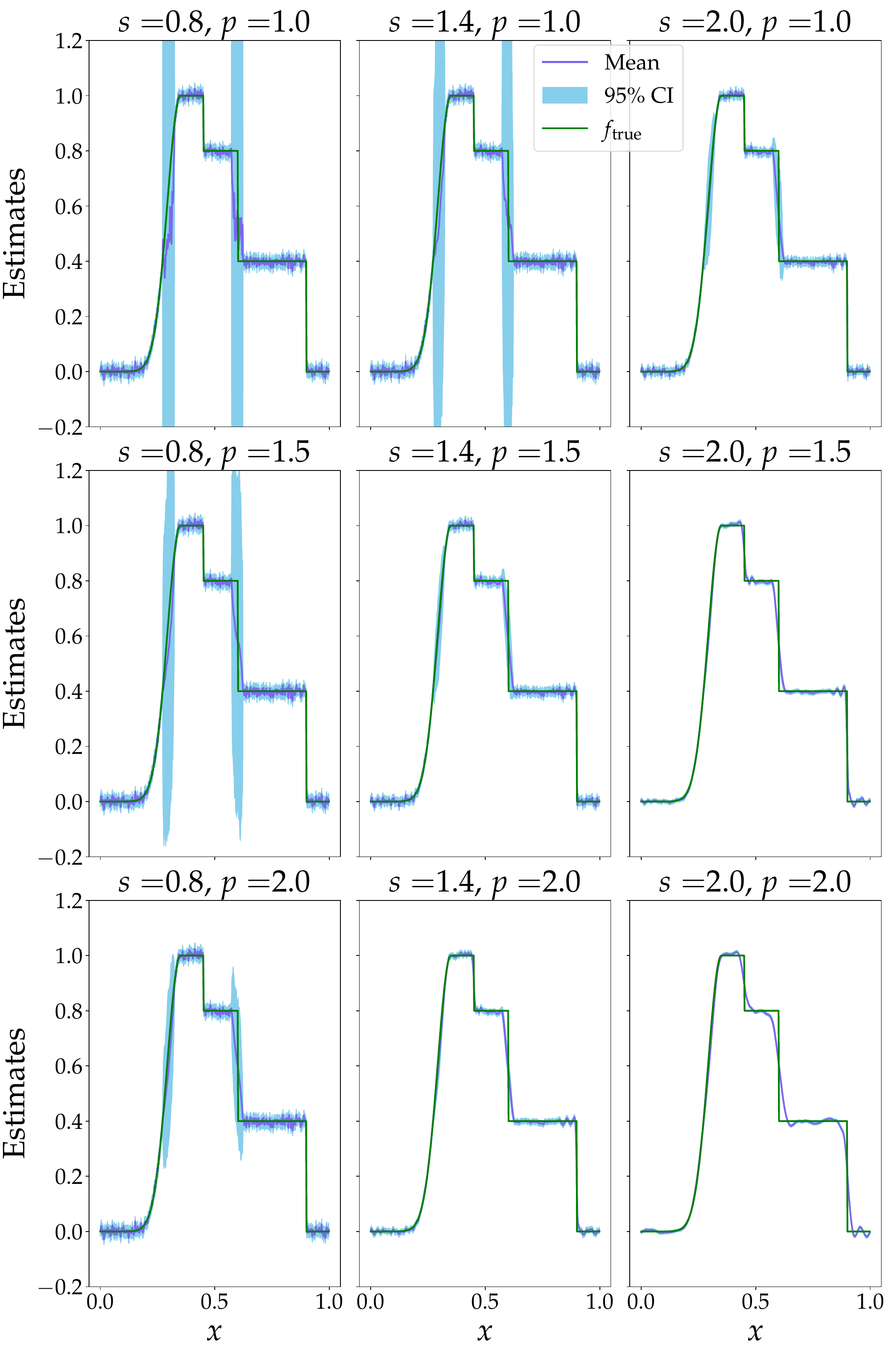}
     \caption{Estimates of a posterior with db8 Besov prior where the parameters $s$ and $p$ are varied between $\{0.8,1.4,2.0\}$ and $\{1.0,1.5,2.0\}$, respectively.}
     \label{fig:ExploringParameters}
 \end{figure}

 The second inpainting test scenario investigates how the choice of Besov parameters $s$ and $p$ affects the estimates of the posterior. In this test, we use the db8 wavelet and set $\lambda=0.025$. We vary the Besov parameters among any combination of $s=\{0.8,1.4,2.0\}$ and $p=\{1.0,1.5,2.0\}$. For each configuration of the parameters, we do the same computations of the posterior estimates as in Section~\ref{sec:ChangeWavelet}.
 The resulting parameter chart of posterior estimates is visualized in Figure~\ref{fig:ExploringParameters}. It is clear that in most cases the uncertainties in the inpainting regions are much larger than in the remaining regions. When we increase $s$ the posterior 95\% CI decreases drastically, and the posterior becomes more concentrated around the mean. Increasing $p$ also decreases the posterior 95\% CI at a similar rate as for increasing $s$. The posterior concentration around the mean with increasing parameters can partially be explained by the faster decay of the deterministic coefficients of the Besov prior in \eqref{eq:BesovExpansion} limiting the variability of the posterior. Additionally, we observe that increasing the parameters $s$ and $p$ results in a smoother posterior mean estimate, which aligns well with the ground truth in the smooth inpainting region. However, the discontinuous part in the inpainting region cannot be well restored, which is due to the choice of the wavelet basis as discussed in the previous subsection. Even though the parameters $s$ and $p$ are of different nature in the definition of the Besov prior we do not observe a significant difference in the results between changing $s$ or $p$.

 \subsection{1D Deconvolution}

 In this section, we consider the 1D deconvolution problem. Under this test problem, we illustrate the discretization invariance property of the Besov prior and compare the efficiency of the RTO-MH sampler with another state-of-the-art sampler, the No-U-Turn sampler (NUTS) Hamiltonian Monte Carlo method \citep{Hoffman2014}. Note that due to the rapid decay of the singular values in the forward operator from the inpainting problem, NUTS is not applicable to the problem in the previous subsection but can be applied here.
 
\begin{figure}[t]
     \centering
     \includegraphics[width=76.0mm]{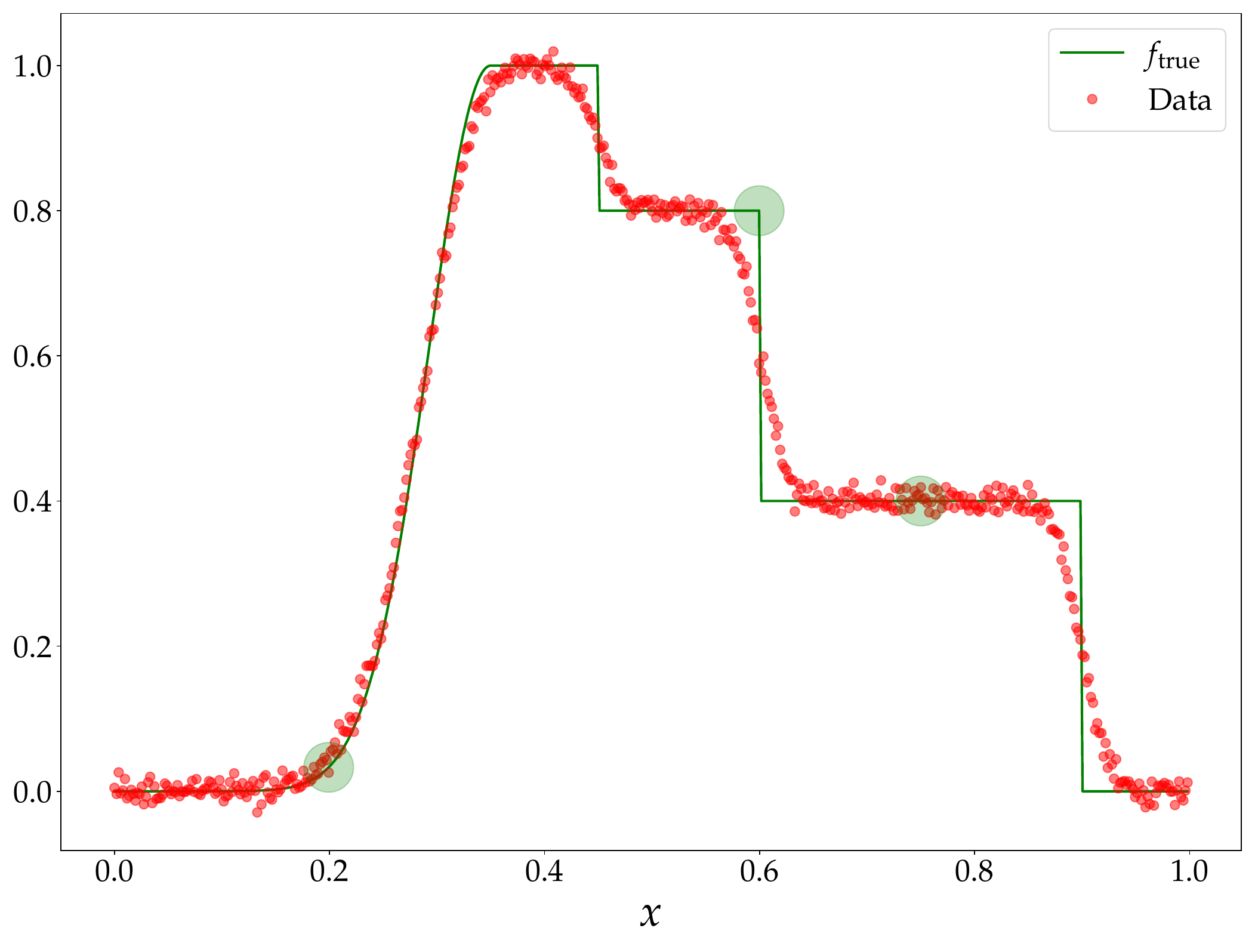}
     \caption{Ground truth signal, noisy convolution data, and the green points are at the positions $x=0.2$, $x=0.6$, and $x=0.75$ which is where we compare the computed chains and autocorrelation functions.}
     \label{fig:Convolutiondata}
 \end{figure}

The periodic deconvolution problem is an inverse problem with forward operator 
 \begin{equation*}
     A(f)(x)=\int_{\mathbb{T}}K(x-y)f(y)dy,
 \end{equation*}
 where $K\in C^{\infty}\left(\mathbb{T}\right)$ is the convolution kernel. In this test, the convolution kernel is chosen to be a 1-periodic Gaussian with zero mean and the standard deviation $\sigma_{\text{ker}}>0$, given by
 \begin{equation*}
      K(x)=\frac{1}{\sqrt{2\pi}\sigma_{\text{ker}}}\sum_{k=-\infty}^{\infty}\exp\left(\frac{-(x+k)^{2}}{2\sigma_{\text{ker}}^{2}}\right).
 \end{equation*}
 In this test, we choose the discrete forward operator $A\in \mathbb{R}^{n\times n}$ in \eqref{eq:DiscreteInverseProblem} to perform discrete periodic convolution using a Gaussian kernel truncated at $\pm 3\sigma_{\text{ker}}$.
 We choose the ground truth signal as in Section~\ref{sec:Inpainting} and set the kernel standard deviation to $\sigma_{\text{ker}}=0.02$. The noise $\epsilon_{n}$ is additive white Gaussian noise with mean $0$ and standard deviation $\sigma>0$, chosen such that $\frac{\|\epsilon_{n}\|_{2}}{\|f_{\text{true}}\|_{2}}=0.02$. 
 The ground truth and the degraded data $y_n$ with $n=512$ are shown in Figure~\ref{fig:Convolutiondata}.

\subsubsection{Discretization invariance}

To illustrate the discretization invariance property from the Besov prior, we discretize the interval $[0,1]$ on the grid $x_{n}=\frac{i}{n}$, $i=0,1,\ldots,n-1$ with $n\in \{32,64,128,256,512,1024,2048\}$.  
In addition, we set the Besov parameters to $\lambda=1.0$, $s=1.0$ and $p=1.5$, and use two Besov prior configurations, i.e., one with the Haar wavelet basis and one with the db8 wavelet basis. 

\begin{figure}[t]
    \centering
    \includegraphics[width=76.0mm]{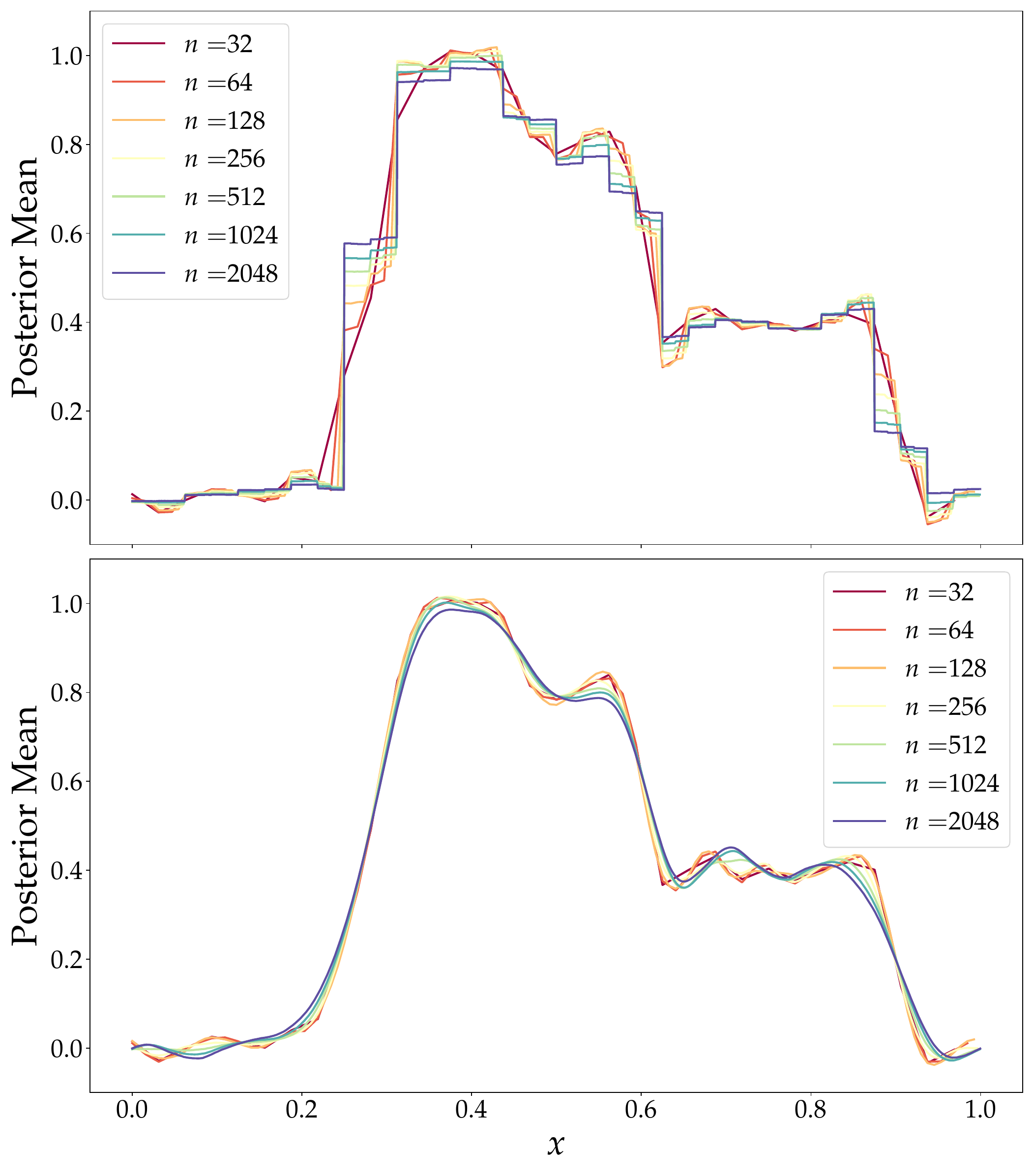}
    \caption{Estimated posterior means with Haar wavelet basis (in the top) and db8 wavelet (in the bottom) at various resolutions $n\in \{32,64,128,256,512,1024,2048\}$.}
    \label{fig:Discretization_Invariance}
\end{figure}

Figure~\ref{fig:Discretization_Invariance} shows the estimates of the posterior mean from 3000 accepted samples in each case, obtained using the RTO-MH sampler proposed in the paper. It is clear that, for both wavelets, the posterior mean converges as we increase the resolution of the signal. This provides numerical evidence for the discretization invariance of the Besov prior as proved in \citep{Lassas2009}.

 \subsubsection{Comparison with NUTS}
In this section, we compare the computational efficiency of the proposed RTO-MH sampler with a state-of-the-art sampler, the  NUTS Hamiltonian Monte Carlo method, where we use the CUQIpy NUTS implementation \citep{Riis2023} for computations. We consider the Haar wavelet Besov prior with parameters $\lambda=1.0$, $s=1.0$ and $p=1.5$ and use both samplers to compute 1000 samples from the resulting posterior. For NUTS we compute an additional 400 burn-in samples to tune the sampler. The RTO-MH sampler does not need burn-in, and we compute 1400 samples to achieve 1000 accepted samples. 

 \begin{figure}[t]
     \centering
     \includegraphics[width=76.0mm]{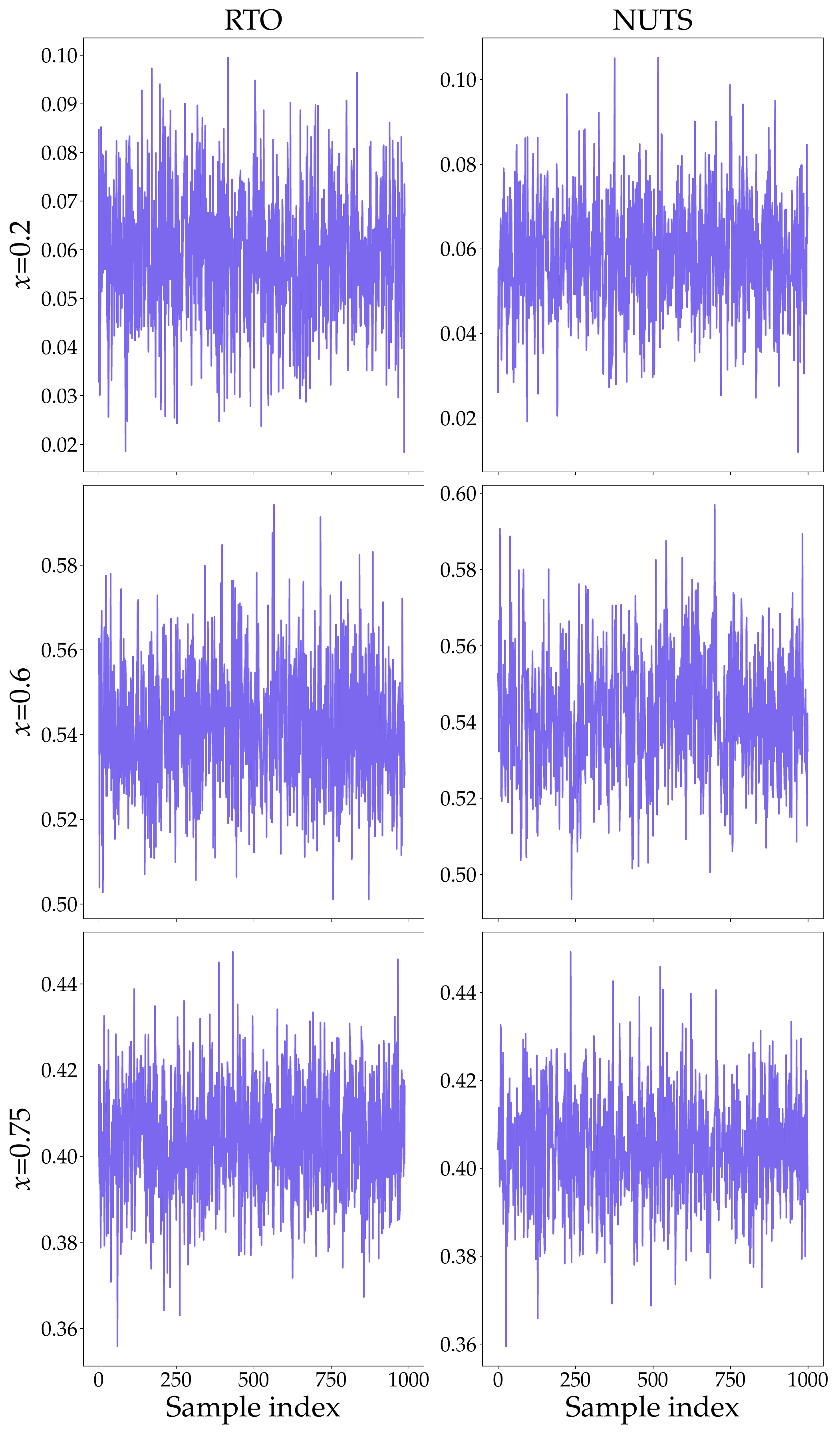}
     \caption{MCMC chains from the points at $x=0.2$, $x=0.6$ and $x=0.75$ with each sampler.}
     \label{fig:Comparison_chains}
 \end{figure}
 \begin{figure}[t]
     \centering
     \includegraphics[width=76.0mm]{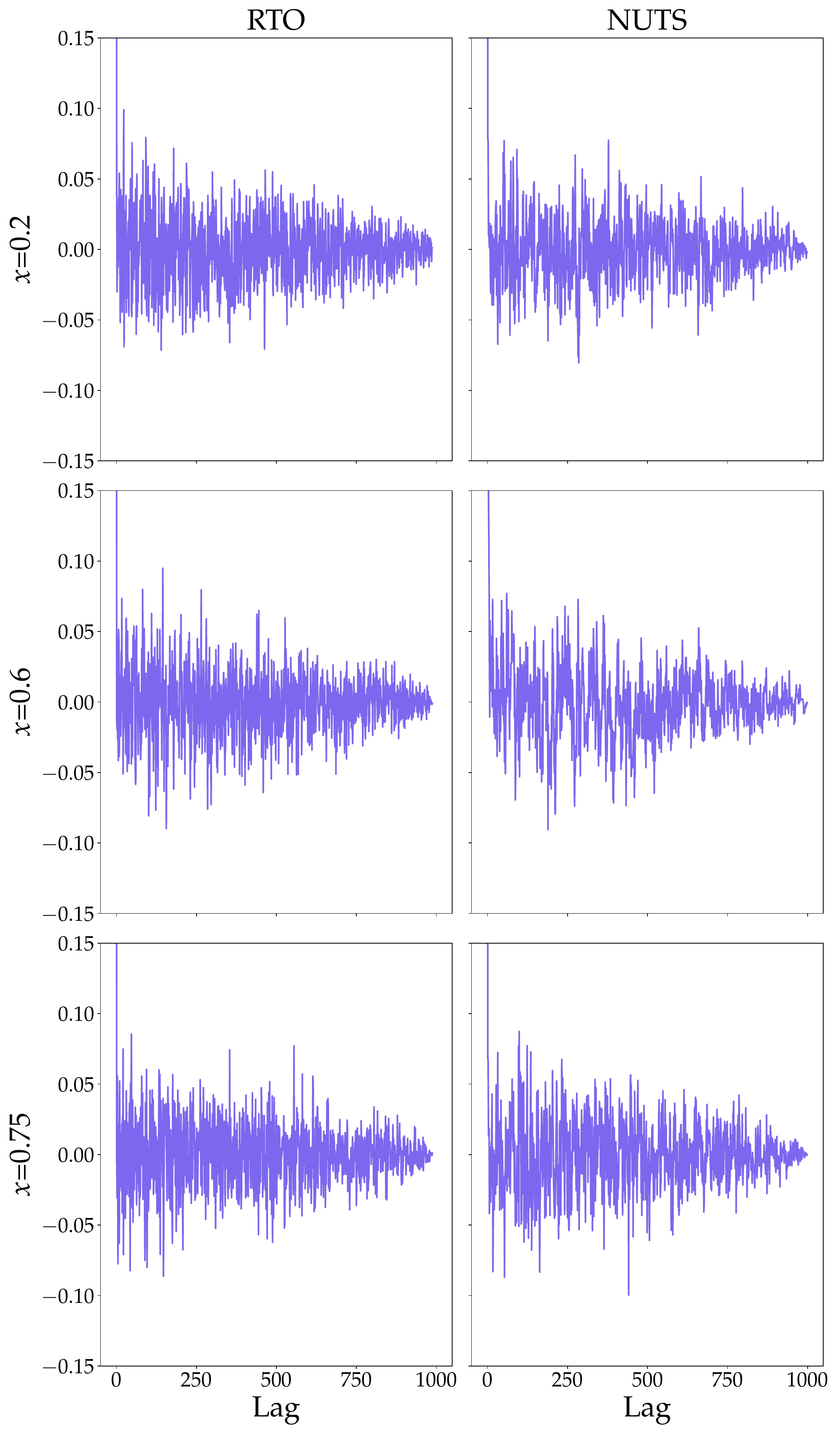}
     \caption{MCMC ACF plots from samples at the points $x=0.2$, $x=0.6$, and $x=0.75$ for each sampler.}
     \label{fig:Comparison_acf}
 \end{figure}

\begin{figure}[t]
    \centering
    \includegraphics[width=76.0mm]{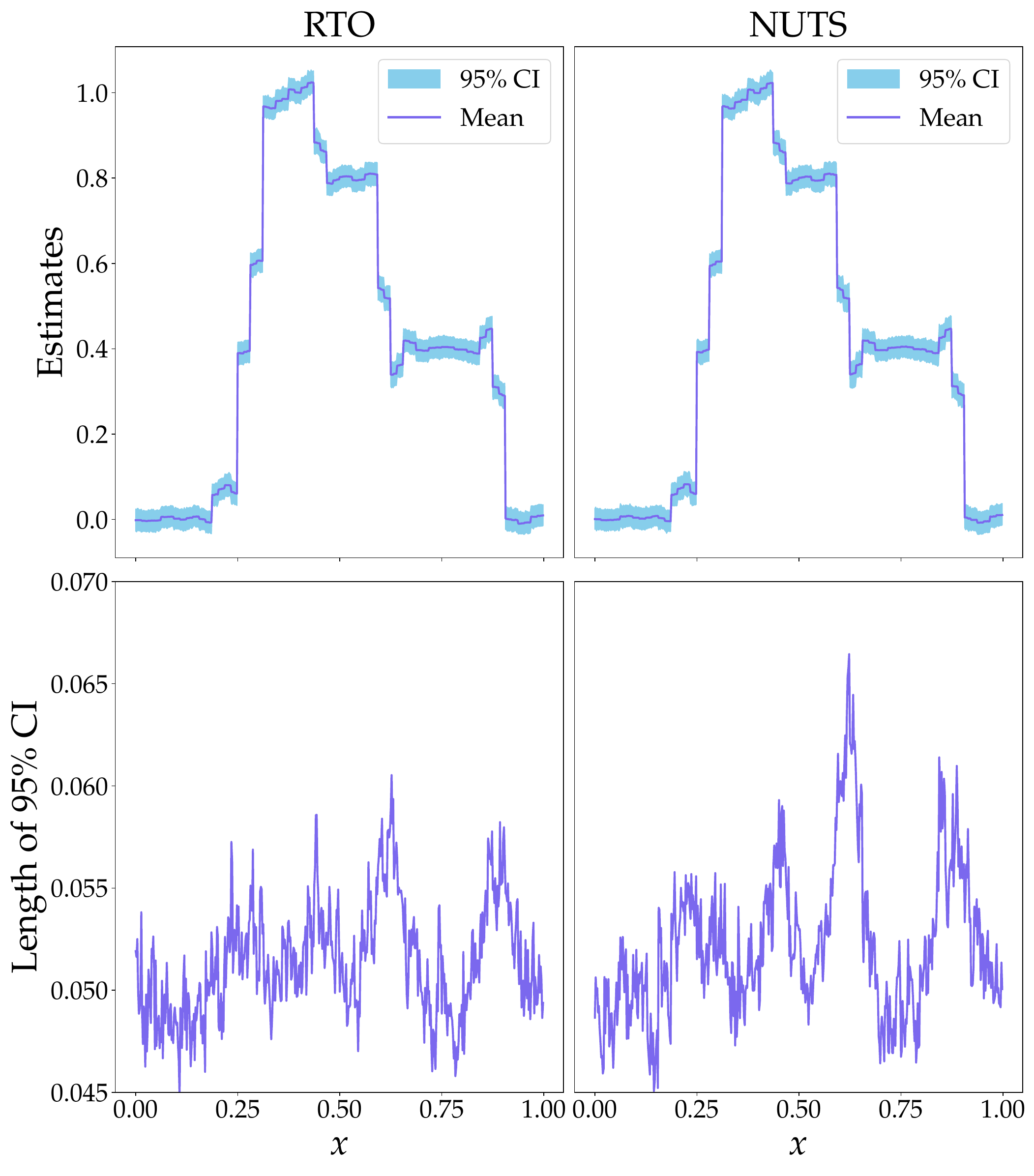}
    \caption{Comparison between the RTO-MH sampler and NUTS with respect to the estimates of the posterior mean and the width of the 95\% CI.}
    \label{fig:Comparison_uq}
\end{figure}

Three computed chains and their autocorrelation function (ACF) are shown in Figure~\ref{fig:Comparison_chains} and Figure~\ref{fig:Comparison_acf}. The chains correspond to the three locations highlighted in Figure~\ref{fig:Convolutiondata}, and they are chosen at three different locations: from a smooth region, a corner and from a constant region, respectively. The Markov chains look stationary and have good mixing which is also shown in the ACF plots that decrease rapidly. In Figure~\ref{fig:Comparison_uq} we show the posterior means and the 95\% CI estimates from both samplers. We can see that the estimated posterior means coincide and the width of the 95\% CI deviates slightly. The slight difference may be due to the number of samples not being large enough to estimate the standard deviation precisely as well as the sample correlation obtained by each sampler. 

We measure the computational efficiency in terms of effective sample size (ESS) per CPU time. ESS is the number of samples used by an independent Monte Carlo estimator with the same variance as the estimator computed by the correlated MCMC samples. We compute the estimated ESS as
 \begin{equation*}
     \hat{N}_{\text{ESS}}=\frac{MN}{-1+2\sum_{t=0}^{K}\hat{\rho}_{2t}+\hat{\rho}_{2t+1}},
 \end{equation*}
 where $M$ is the number of chains, $N$ is the amount of samples, $\hat{\rho}_{t}$ is the estimated autocorrelation function and $K$ is the last integer for which $\hat{\rho}_{2K}+\hat{\rho}_{2K+1}>0$. We compute the ESS and ESS per CPU-second for every MCMC chain and tabulate the minimum, median, and maximum in Table~\ref{tab:ESS_comparison}. We can see that the RTO-MH sampler has higher ESS compared with NUTS. According to the ESS/CPU-second values, the RTO-MH sampler and NUTS are rather identical concerning computational efficiency.  
\begin{table}[t]
\resizebox{76.0mm}{!}{
\begin{tabular}{|c|ccc|ccc|}
\hline
Method & \multicolumn{3}{c|}{ESS}    & \multicolumn{3}{c|}{ESS /CPU-second} \\ \hline
       & Minimum & Median  & Maximum & Minimum    & Median    & Maximum   \\ \hline
RTO-MH & 580.65  & 954.15 & 1178.22 & 0.77       & 1.26      &  1.56      \\
NUTS   & 240.21  & 677.52 & 1064.86 & 0.38       &  1.07      & 1.68     \\ \hline
\end{tabular}
}
\caption{ESS and ESS/CPU-second comparison between RTO-MH and NUTS.}
\label{tab:ESS_comparison}
\end{table}

 \subsection{2D Computed Tomography}
 
 Computed Tomography (CT) is an imaging technique that reconstructs cross-sections or slices of an unknown object from external measurements of intensity attenuation of X-rays. The unknown, often referred to as the attenuation field, is usually assumed to be piecewise constant (e.g. in medical imaging) reflecting groups of different materials with constant attenuation coefficients. The most popularly used reconstruction method for CT is the filtered back projection (FBP) \citep{hansen2021computed}. However, FBP tends to produce reconstruction artifacts when the data are limited in the sense of sparse measurements or limited projection angles. Total variation (TV) regularization is another widely used method for promoting piecewise constant features in CT and can handle limited data settings. However, TV does not fit into the Bayesian framework, which poses challenges in interpreting the uncertainties \citep{Matti_Lassas_2004}.
 As an alternative, we suggest using the Besov priors to promote piecewise constant features that are consistent with the Bayesian framework.
 
 The interaction of X-ray with the attenuation field can be modelled using line integrals. The transformation of an attenuation field $f$ into its line integrals along straight lines is given by a linear transformation known as the Radon-transform \citep{Radon1986}. 
 The Radon-transform of a 2D function $f$ is given by
 \begin{align*}
     \mathcal{R}f(\theta,\eta)&=\int_{L(\theta,\eta)}f(x(t),y(t))dL(t) \nonumber\\
     &=\int_{-\infty}^{\infty}f(t\sin(\theta)+\eta\cos(\theta),\nonumber\\&\qquad \qquad \;-t\cos(\theta)+\eta\sin(\theta))dt,
 \end{align*}
 where $\theta\in [0, \pi[$ is the angle between the normal vector of the line $L(\theta,\eta)$ and the $x$-axis and $\eta\in [-1,1]$ is the signed distance between $L(\theta,\eta)$ and the origin. 

 The discrete 2D Radon-transform is performed for a finite set of lines in the plane. 
The parallel-beam geometry, used in this section, comprises a uniform selection of the parameters $\theta$ and $\eta$. 
 
\begin{figure}[t]
     \centering
     \begin{subfigure}[b]{76.0mm}
         \centering
         \includegraphics[width=\textwidth]{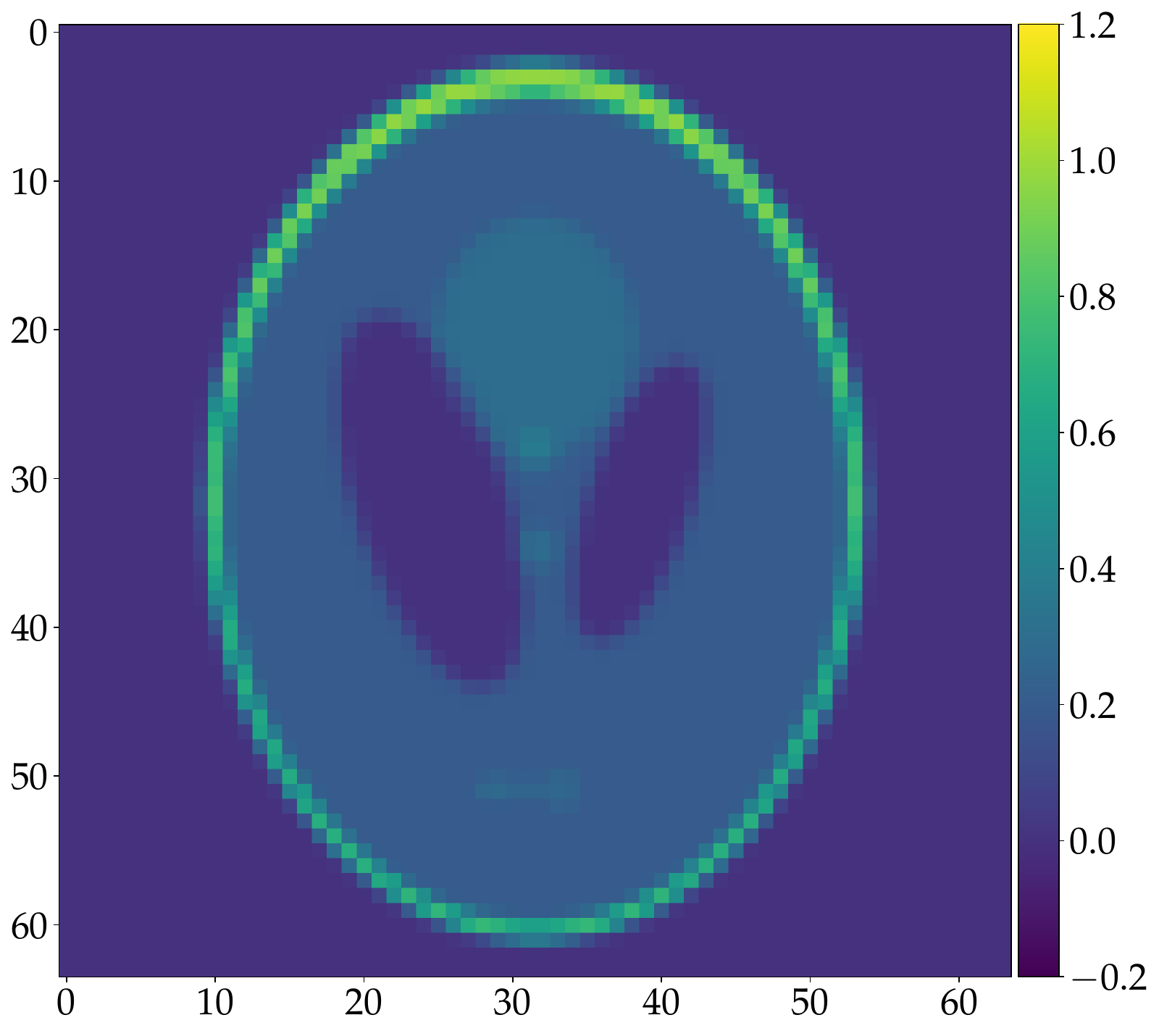}
         \caption{Sheep-Logan}
     \end{subfigure}
     \begin{subfigure}[b]{76.0mm}
         \centering
         \includegraphics[width=0.475\textwidth]{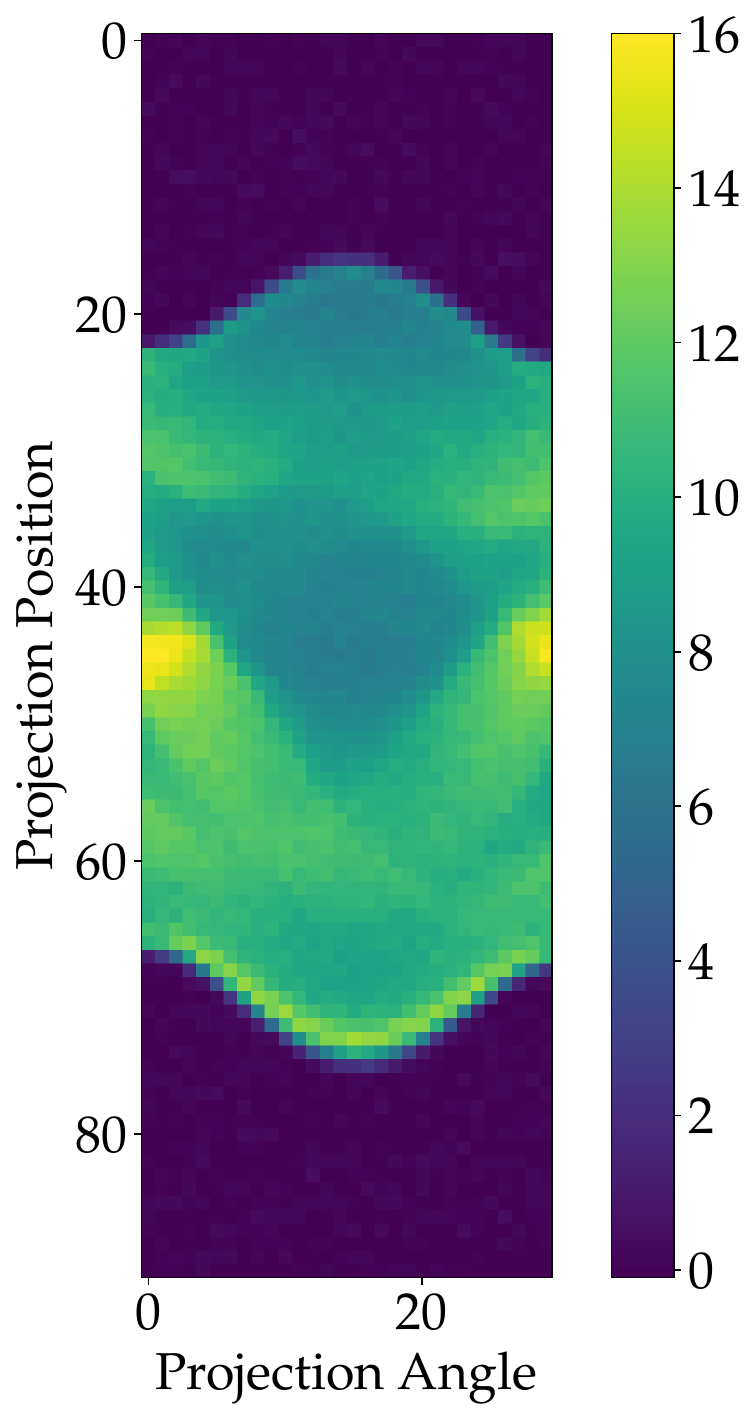}
         \caption{CT data}
     \end{subfigure}
        \caption{Ground truth of the CT problem together with the corresponding noisy data.}
        \label{fig:CTdata}
\end{figure}

 We denote $A$ to be the discrete Radon-transform as the forward operator for the CT problem in \eqref{eq:DiscreteInverseProblem}. The output of $A$ is referred to as a sinogram. We impose a Besov prior to the CT problem by assuming that the continuous unknown attenuation field $f$ follows a Besov prior as in \eqref{eq:FiniteBesovExpansion}. We refer to the discrete attenuation field as $f_{n}$. To set up the test problem, we discretize the tangent $\theta \in [0,\pi[$ on an equidistant grid of $30$ points and the intercept $\eta\in [-1,1]$ on an equidistant grid of $91$ points, respectively. Therefore, the dimension of the sinogram is $m=30\times 91$ and this represents a sparse data scenario known as sparse angle CT. We use the Python package scikit-image \citep{scikit-image} to compute the discrete Radon-transform. We choose the true unknown attenuation field to be the Sheep-Logan phantom with the dimension $64\times 64$, shown in Figure~\ref{fig:CTdata}, accessed from the scikit-image package. The relative noise level is set as $0.02$ as defined in Section~\ref{sec:Inpainting}.

\begin{figure}[t]
     \centering
     \begin{subfigure}[b]{76.0mm}
         \centering
         \includegraphics[width=76.0mm]{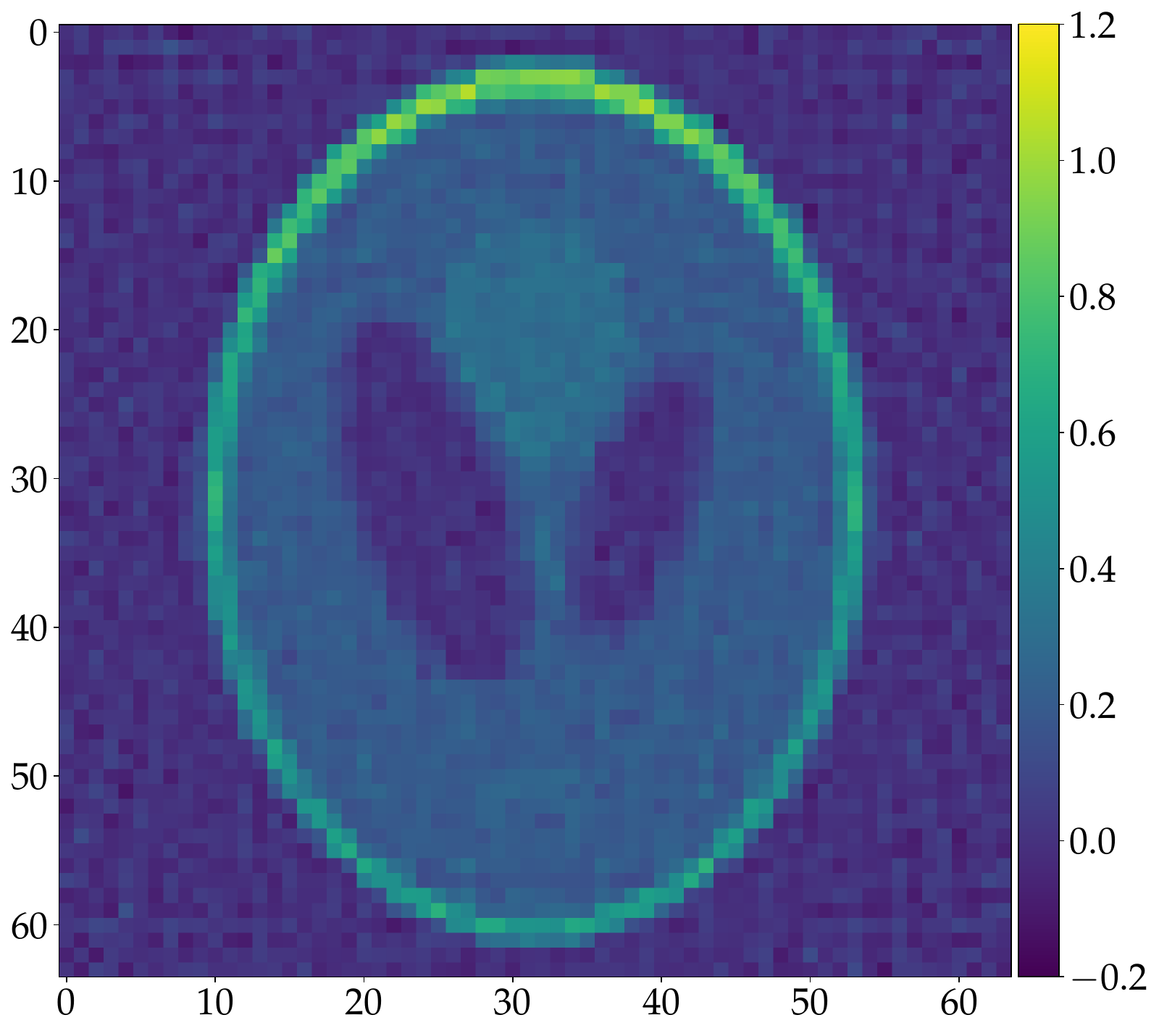}
         \caption{Posterior Mean}
     \end{subfigure}
     \hfill
     \begin{subfigure}[b]{76.0mm}
         \centering
         \includegraphics[width=76.0mm]{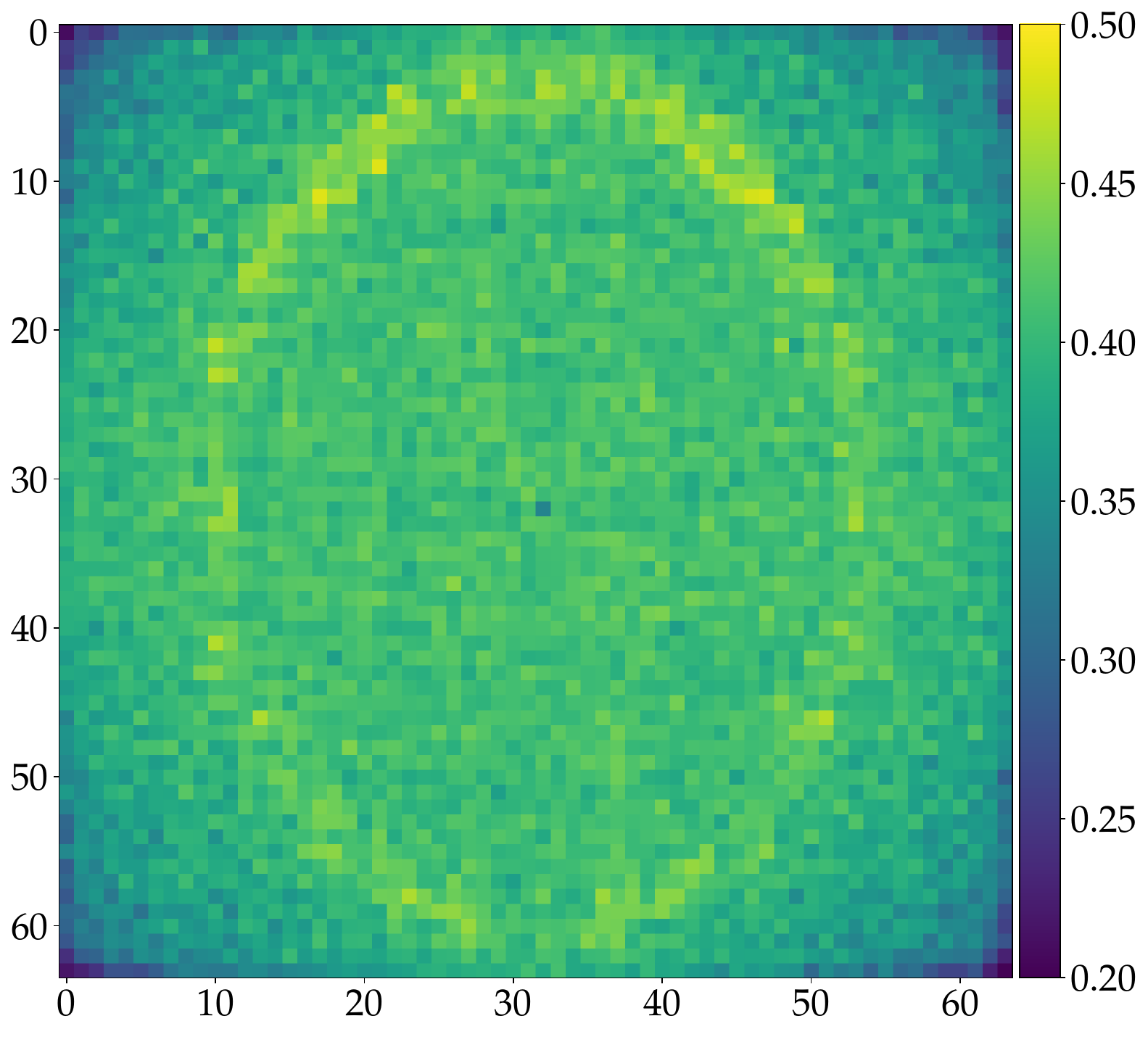}
         \caption{Pixelwise width of 95\% posterior CI}
     \end{subfigure}
        \caption{The estimated Posterior mean and the estimated pixelwise width of the 95\% posterior credible interval.}
        \label{fig:CTposterior}
\end{figure}

We choose a Besov prior with Haar wavelet basis and parameters $s=1.0$, $p=1.5$ and $\lambda=0.025$. This choice promotes piecewise constant samples and reveals the uncertainty in the estimate. We apply the proposed RTO-MH sampler to draw $5\times 10^4$ samples from the posterior defined in \eqref{eq:DiscretePosterior}, of which 2315 samples were accepted.  
Fig.~\ref{fig:CTposterior} shows the computed posterior mean and the pixelwise width of the 95\% credible interval. We see that the posterior mean can reconstruct the piecewise constant phantom, although it has small noisy artifacts. Furthermore, we notice that the 95\% posterior credible interval is most significant around the highest discontinuities while we do not see a general pattern on the rest of the domain.

\begin{figure}[t]
    \centering
    \includegraphics[width=76.0mm]{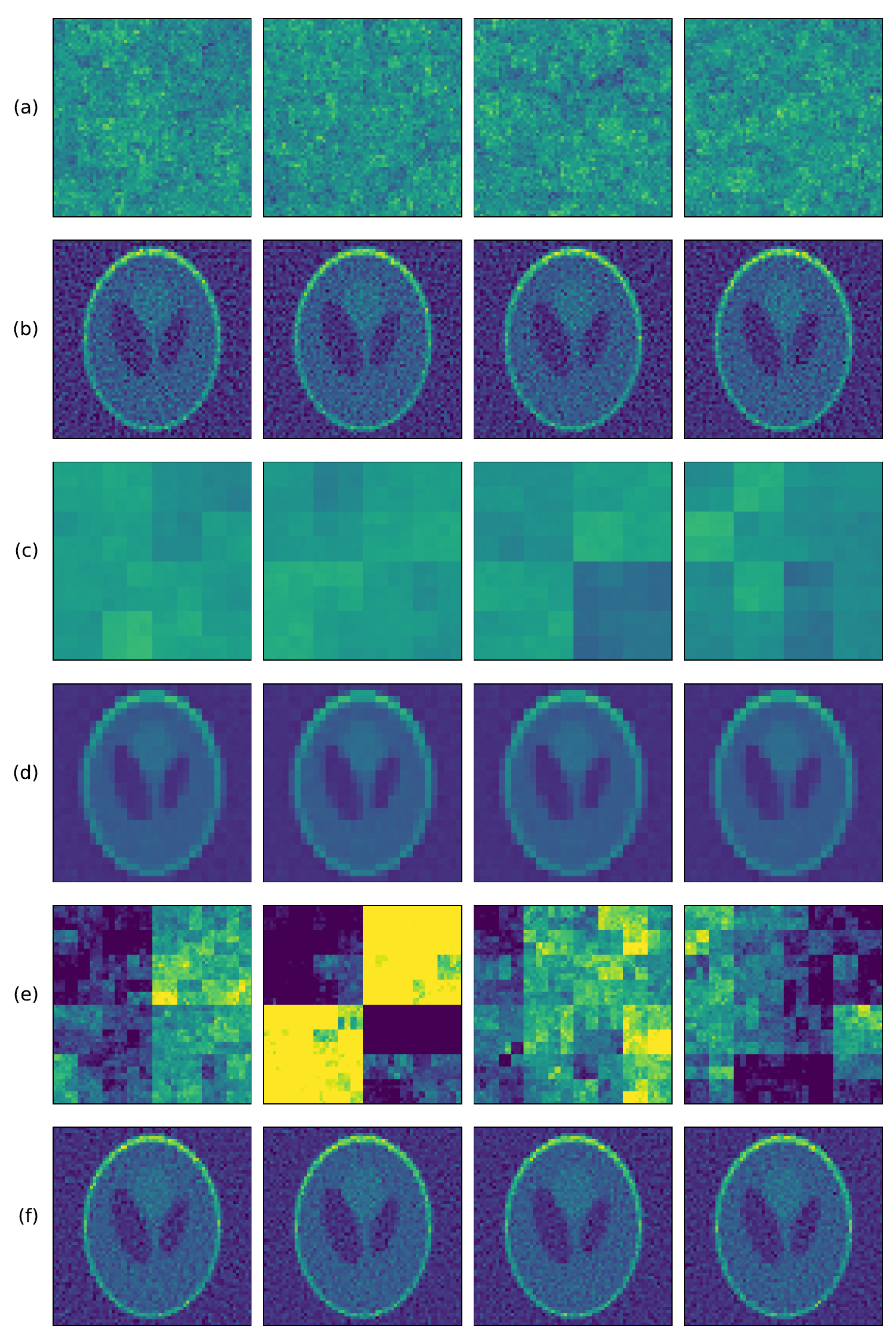}
    \caption{(a) and (b) are prior samples and posterior samples in the case where $s=1.0$, $p=1.5$ and $\lambda=0.025$, respectively. (c) and (d) are prior samples and posterior samples in the case where $s=2.5$, $p=1.5$ and $\lambda=0.025$, respectively. (e) and (f) are prior samples and posterior samples in the case where $s=2.5$, $p=1.0$ and $\lambda=0.025$, respectively.}
    \label{fig:prior_posterior_samples}
\end{figure}

We also provide prior and posterior samples for three cases of prior parameters $(s=1.0,p=1.5),(s=2.5,p=1.5),$ and $(s=2.5,p=1.0)$ with the same $\lambda$ in Fig.~\ref{fig:prior_posterior_samples}. The parameter $\lambda$ was chosen such that all the visualized prior samples are descriptive of the unknown, that is, the interval width of pixel values of the prior samples is at least 0.5. It is clear that increasing $s$ and $p$ results in smoother prior as well as posterior samples. Further, by comparing the prior samples, we can see the distinct influence of changing $s$ and $p$. With the same $p$ by increasing $s$, i.e. comparing rows (a) and (c) in Fig.~\ref{fig:prior_posterior_samples}, the prior samples are mainly represented by large smoothing features, which is due to the significant decreasing of the wavelet coefficients according to small features. Comparing rows (e) and (c) from Fig.~\ref{fig:prior_posterior_samples}, we can see that increasing $p$ results in prior samples with a globally decreased range across all levels of features.

\begin{figure}[t]
     \centering
     \includegraphics[width=76.0mm]{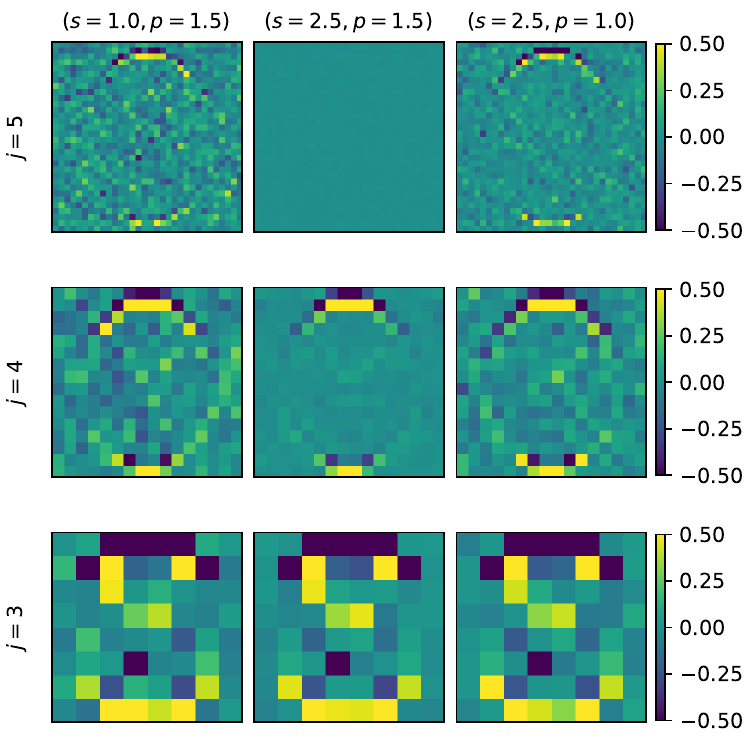}
     \caption{Horizontal $(\ell=1)$ wavelet coefficients of a posterior sample for all three choices of parameters. In each row $j$ represents the wavelet level as in \eqref{eq:FiniteBesovExpansion}, and columns represent the choice of prior parameters.}
     \label{fig:Coeff_1}
 \end{figure}

\begin{figure}[t]
     \centering
     \includegraphics[width=76.0mm]{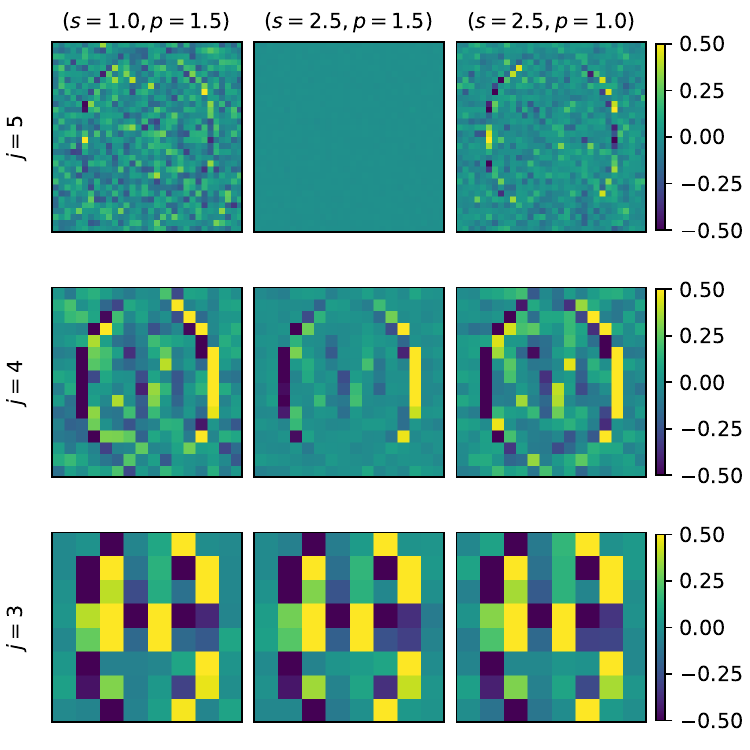}
     \caption{Vertical $(\ell=2)$ wavelet coefficients of a posterior sample for all three choices of parameters. In each row $j$ represents the wavelet level as in \eqref{eq:FiniteBesovExpansion}, and columns represent the choice of prior parameters.}
     \label{fig:Coeff_2}
 \end{figure}

To further illustrate the different influence from changing $s$ and $p$, we provide the value of the wavelet coefficients from different scaling levels $j$ with the horizontal and vertical wavelets
for a single posterior sample in Fig.~\ref{fig:Coeff_1} and~\ref{fig:Coeff_2}. It is evident that the coefficients become sparser at higher levels, and increasing $s$ results in even sparser coefficients at these levels, as $s$ governs the decay of the wavelet coefficients. In addition, increasing $p$ promotes higher sparsity across all levels rather uniformly.
We report that this pattern is repeated in other posterior samples. 


%\section{Discussion}\label{sec5}

\section{Conclusion}\label{sec5}
This paper presents a Bayesian framework for reconstructing piecewise constant functions by using Besov priors in linear inverse problems. To enable inference within this framework, we addressed the issue of sampling from a difficult and non-Gaussian posterior, using an RTO approach. Through a series of numerical experiments, we demonstrated that the RTO method is comparable to the NUTS method in efficiently sampling a posterior and can be applied to more types of inverse problems such as image inpainting. In addition, we show that our method is discretization invariant, which makes it suitable for exploring high-dimensional posteriors. Beyond sampling efficiency, our study also highlights desirable properties of Besov priors, including their ability to promote smoothness depending on the choice of wavelet basis. Our results further suggest that Besov priors promote sparsity in the posterior samples depending on the choices of prior parameters.  

While the current work focuses on synthetic examples, designed to showcase key features of the Besov prior and the proposed RTO sampling method, it does not yet address the issue of scalability. The main computational bottleneck of the RTO method is the cubic scaling $\mathcal{O}(n^{3})$ of computing the absolute value of the determinant of $Q^{T}J_{\mathcal{A}}(h_{n})$ in eq. \eqref{eq:MHweights} for each sample. Consequently, the RTO method cannot be directly applied to a high-dimensional real-sized inverse problem. Related work, such as \citep{Wang2020}, addresses the scalability limitations by using techniques from numerical linear algebra. Building on these advances, a natural step is to explore whether such techniques combined with parallelization can make our RTO method applicable to high-dimensional real-world inverse problems.

The application of MCMC methods to high-dimensional inverse problems is commonly done using algorithms with low per-sample cost, e.g., preconditioned Crank-Nicolson \citep{cotter2013} or Langevin based methods \citep{Durmus2019}. For an in depth comparison between RTO and Langevin based algorithm we refer the reader to \citep{Laumont2025}.

Interesting directions of future research would be to consider the prior parameters $s$ and $p$ as hyperparameters resulting in a hierarchical structure. Another promising direction would be to extend our framework to suitable nonlinear inverse problems and study the usefulness of Besov priors in this context. 

%\backmatter

%\bmhead{Supplementary information}

%\bmhead{Acknowledgments}

%\section*{Declarations}

%\begin{appendices}

%\section{Appendix}\label{secA1}




%\end{appendices}

\bibliography{Bibliography}% common bib file

\end{document}